\documentclass[12pt,amsymb,fullpage]{amsart}
\usepackage{amssymb,amscd,pstricks}

\newtheorem{theorem}{Theorem}[section]
\newtheorem{defn}[theorem]{Definition}

\newtheorem{lemma}[theorem]{Lemma}
\newtheorem{fact}[theorem]{Fact}

\newtheorem{eple}[theorem]{Example}
\newtheorem{rmk}[theorem]{Remarks}
\newtheorem{dsc}[theorem]{Discussion}
\newtheorem{nota}[theorem]{Notation}

\newsavebox{\indbin}
\savebox{\indbin}{\begin{picture}(0,0)
\newlength{\gnu}
\settowidth{\gnu}{$\smile$} \setlength{\unitlength}{.5\gnu}
\put(-1,-.65){$\smile$} \put(-.25,.1){$|$}
\end{picture}}

\newcommand{\be}{\begin{enumerate}}
\newcommand{\bd}{\begin{defn}}
\newcommand{\bt}{\begin{theorem}}
\newcommand{\bl}{\begin{lemma}}
\newcommand{\ee}{\end{enumerate}}
\newcommand{\ed}{\end{defn}}
\newcommand{\et}{\end{theorem}}
\newcommand{\el}{\end{lemma}}

\begin{document}
\title{A Theory of Divisors for Algebraic Curves}
\author{Tristram de Piro}
\address{55b, Via Ludovico Albertoni, Rome, Italy}
\thanks{Thanks to Francesco Severi and The Lamb}
\begin{abstract}
The purpose of this paper is two-fold. We first prove a series of results, concerned with the notion of Zariski multiplicity, mainly for non-singular algebraic curves. These results are required in \cite{depiro1}, where, following Severi, we introduced the notion of the "branch" of an algebraic curve. Secondly, we use results from \cite{depiro1}, in order to develop a refined theory of $g_{n}^{r}$ on an algebraic curve. This refinement depends critically on replacing the notion of a point with that of a "branch". We are then able to construct a theory of divisors, \emph{generalising} the corresponding theory in the special case when the algebraic curve is \emph{non-singular}, which is \emph{birationally invariant}.
\end{abstract}

\maketitle

\begin{section}{Introduction}
In this paper, we use the same definition of an algebraic curve as in \cite{depiro1}. Namely, an algebraic curve $C$ is a closed, irreducible subvariety of dimension $1$ in $P^{w}$, for some $w\geq 1$, where $P^{w}$ denotes projective space of dimension $w$. We will often abbreviate the terminology of "algebraic curve" to just "curve". The advantage of developing a birationally invariant theory of divisors for such curves depends mainly on the viewpoint of the "Italian School" of algebraic geometry. Namely, that there are a number of benefits in studying the geometry of \emph{plane} algebraic curves, $(*)$, and that any algebraic curve $C$ is birational to a plane algebraic curve $C'$, (see, for example, Theorem 1.33 of \cite{depiro1}). It is not the purpose of this paper to discuss the question raised in $(*)$, leaving this point of view for another occasion. The results of this paper cover \emph{all} characteristics of the underlying algebraically closed field $L$. However, we will make it clear when a result depends on the assumption that $L$ has non-zero characteristic.

\end{section}

\begin{section}{Smooth Curves}

Before looking at this section, the reader is strongly advised to consult the paper \cite{depiro2} for relevant notation and terminology. In particular, the reader should be acquainted with the statement of Theorem 3.3 from \cite{depiro2}. We first recall the following theorem (which was Theorem 6.5 in \cite{depiro2});\\

\begin{theorem}

Let hypotheses be as in Theorem 3.3 of \cite{depiro2}, with the additional
assumption that $char(L)=0$ and $F$, $D$ are \emph{smooth} curves. Then the
notions of Zariski multiplicity and algebraic multiplicity coincide.

\end{theorem}

As we need to refer to the proof of this result from \cite{depiro2} later in the paper, for the convenience of the reader, we repeat it below. The reader should, however, consult \cite{depiro2} for relevant notation.\\

\begin{proof}

As $D$ has a non-constant meromorphic function, we can write $D$ as a finite cover of $P^{1}(L)$.  As we have
 checked both algebraic multiplicity and Zariski multiplicity are multiplicative over composition (in \cite{depiro2}), a straightforward calculation shows that we need only check the notions agree for the branched finite cover $\pi:F\rightarrow P^{1}(L)$. (1)\\

Now consider this cover restricted to $A^{1}$, let $x$ be the canonical coordinate  with $ord_{a} (\pi^{*}(x))=m$, so we have that $\pi^{*}x=h^{m}u$ , for $u$ a unit in ${\mathcal{O}}_{a}$ and $h$ a uniformiser at $a$. (2)\\

As $u$ is a unit and $char(L)=0$, the equation $z^{m}=u$ splits in the residue field of ${\mathcal O}^{\wedge}_{a}$. By
 Hensel's Lemma and Theorem 5.5 of \cite{depiro2}, it is solvable in ${\mathcal O}_{a}^{\wedge}$. By the definition of
 ${\mathcal{O}}_{a}^{\wedge}$, we can find an etale morphism $\pi:(U,b)\rightarrow (F,a)$ containing such a solution
 in the local ring ${\mathcal{O}}_{b}$. We may assume that $U$ is irreducible and moreover, as $\pi$ is etale, that $U$
 is smooth. (3)\\

Now we can embed $U$ in a projective smooth curve $F'$ and, as $F'$ is smooth,  extend the morphism $\pi$ to a projective morphism from $F'$ to $F$. (4)\\

We claim that $(ba)\in graph(\pi)\subset F'\times F $ is unramified in the sense of Zariski structures. For this we need the following fact whose algebraic proof relies on the fact that etale morphisms are flat, see \cite{depiro2};\\

\begin{fact}

Any etale morphism can be locally presented  in the form \\

\begin{eqnarray*}
\begin{CD}
V@>g>>Spec((A[T]/f(T))_{d})\\
@VV\pi V  @VV\pi' V\\
U@>h>>Spec(A)\\
\end{CD}
\end{eqnarray*}

where $f(T)$ is a monic polynomial in $A[T]$, $f'(T)$ is invertible in $(A[T]/f(T))_{d}$ and $g,h$ are isomorphisms.   (5)\\
\end{fact}

Using Lemma 4.6 of \cite{depiro2} and the fact that the open set $V$ is smooth, we may safely replace $graph(\pi)$ by
 $\overline {graph (\pi')}\subset F''\times F$ where $F''$ is the projective closure of $Spec((A[T]/f(T))$,
 $F$ is the projective closure of $Spec(A)$ and $\overline {graph(\pi')}$ is the projective closure of $graph(\pi')$
 and show that $(g(b)a)$ is Zariski unramified. Note that over the open subset $U=Spec(A)\subset F$,
  $\overline{graph(\pi')}=Spec((A[T]/f(T)$ as this is closed in $U\times F''$.  For ease of notation, we replace
   $(g(b)a)$ by $(ba)$. (6)\\

Suppose that $f$ has degree $n$. Let $\sigma_{1}\ldots \sigma_{n}$ be the elementary symmetric functions in $n$
 variables $T_{1},\ldots T_{n}$. Consider the equations\\

$\sigma_{1}(T_{1},\ldots, T_{n})=a_{1}$\\

$\ldots$\\

$\sigma_{n}(T_{1},\ldots,T_{n})=a_{n}$ (*)\\

where $a_{1},\ldots a_{n}$ are the coefficients of $f$ with
appropriate sign. These cut out a closed subscheme $C\subset
Spec(A[T_{1}\ldots T_{N}])$. Suppose $(ba)\in
graph(\pi')=Spec(A[T]/f(T))$ is ramified in the sense of Zariski
structures, then I can find $(a'b_{1}b_{2})\in {\mathcal V}_{abb}$
with $(a'b_{1})$,$(a'b_{2})\in Spec(A(T)/f(T))$ and $b_{1},b_{2}$
distinct. Then complete $(b_{1}b_{2})$ to an $n$-tuple
$(b_{1}b_{2}c_{1}'\ldots c_{n-2}')$ corresponding to the roots of
$f$ over $a'$. The tuple $(a'b_{1}b_{2}c_{1}'\ldots c_{n-2}')$
satisfies $C$, hence so does the specialisation $(abbc_{1}\ldots
c_{n-2})$. Then the tuple $(bbc_{1}\ldots c_{n-2})$ satisfies
$(*)$ with the coefficients evaluated at $a$. However such a
solution is unique up to permutation and corresponds to the roots
of $f$ over $a$. This shows that $f$ has a double root at $(ab)$
and therefore $f'(T)|_{ab}=0$. As $(ab)$ lies inside
$Spec(A[T]/f(T))_{d}$, this contradicts the fact that $f'$ is
invertible in $A[T]/f(T))_{d}$. (7)\\

In $(2)$ we may therefore assume that $\pi^{*}x=h^{m}$ for $h$ a local uniformiser at $a$. Now we have the sequence
 of ring inclusions given by \\

$L[x]\rightarrow L[x,y]/(y^{m}-x)\rightarrow R$\\

\ \ \ \ \ \ \ \ \ $x\mapsto \pi^{*}x, y\mapsto h$\\

where $R$ is the coordinate ring of $F$ in some affine neighborhood of $a$. It follows that we can factor our original
 map such that $F$ is etale near $a$ over the projective closure of $y^{m}-x=0$. (8)\\

Again, repeating the argument from (4) to (7), we just need to
check that the projective closure of $y^{m}-x$ has multiplicity
$m$ at $0$ considered as a cover of $P^{1}(\bar k)$. This is
trival, let $\epsilon\in {\mathcal V}_{0}$ be generic over

$\mathcal M$,then as we are working in characteristic $0$ we can
find distinct $\epsilon_{1},\ldots \epsilon_{m}$ in ${\mathcal
M}_{*}$ solving $y^{m}=\epsilon$. By specialisation, each
$\epsilon_{i}\in{\mathcal V}_{0}$. (9)
\end{proof}

The purpose of this section is essentially to find an analogous result to Theorem 2.1 when $char(L)=p\neq 0$. An analogous result was given in \cite{depiro2}, however, the proof was flawed. We correct this difficulty here. We obtained similar results, in \cite{depiro2}, under different assumptions, by the straightforward method of counting points in the fibres. In this section, we need to use more sophisticated local methods, which will be explained below. We first make the following remark concerning the Frobenius morphism;\\

\begin{rmk}{Frobenius}\\

Given a smooth curve $C$, defined over a field of characteristic $p$, with function
field $L(C)$, we let $L(C)^{1/p}$ be the field obtained by
extracting $p^{th}$ roots of $L(C)$ in some fixed algebraic
closure. We denote by $C_{p}$ the unique (up to isomorphism)
smooth curve, having function field $L(C)^{1/p}$. Corresponding to the
inclusion $i:L(C)\rightarrow L(C)^{1/p}$, we obtain a morphism
$Frob:C_{p}\rightarrow C$, which, by some abuse of the standard
terminology, (the standard terminology is $L$-linear Frobenius),
we will refer to as Frobenius. Although $L(C)$ and $L(C)^{1/p}$
are clearly isomorphic as fields, they may not be isomorphic over
$L$. Hence, $C$ and $C_{p}$ are not necessarily isomorphic curves.
The Frobenius morphism may be explicitly realised as follows;\\

Let $C$ be embedded in $P^{n}$, for some $n$, defined by the
homogeneous polynomials $\{f_{1},\ldots,f_{m}\}$. Let $C'$ be the
variety defined by $\{\overline{f_{1}},\ldots,\overline{f_{m}}\}$,
where, for $1\leq j\leq m$, $\overline{f_{j}}$ is the homogeneous
polynomial obtained by applying inverse Frobenius to the
coefficients. Then, by a straightforward calculation using
Jacobians, $C'$ defines a smooth curve. The morphism Frobenius;\\

$Fr:P^{n}\rightarrow P^{n}$\\

$Fr([X_{0}:\ldots:X_{n}])=[X_{0}^{p}:\ldots:X_{n}^{p}]$\\

restricts to define a morphism $Fr:C'\rightarrow C$. Let $Rat_{k}$
denote the rational functions of degree $k$ on $P^{n}$. Then $Fr$
induces a map;\\

$Fr^{*}:Rat_{k}\rightarrow Rat_{kp}$\\

by the formula;\\

$(Fr^{*}F)(X_{0},\ldots,X_{n})=F(X_{0}^{p},\ldots,X_{n}^{p})$\\

For a homogeneous polynomial $f_{j}$ defining $C$, we have that;\\

$Fr^{*}(f_{j})=(\overline f_{j})^{p}$\\

Hence, $Fr^{*}$ restricts to define an $L$-linear map;\\

$Fr^{*}:L(C)\rightarrow L(C')$\\

One can also define a map;\\

$Fr^{-1*}:L[X_{0},\ldots,X_{n}]\rightarrow
L[X_{0}^{1/p},\ldots,X_{n}^{1/p}]$\\

by the formula;\\

$(Fr^{-1*}F)(X_{0},\ldots,X_{n})=F(X_{0}^{1/p},\ldots,X_{n}^{1/p})$\\

For a homogeneous polynomial $\overline{f_{j}}$ defining $C'$, we have that;\\

$Fr^{-1*}(\overline{f_{j}})=(f_{j})^{1/p}$\\

Hence, $Fr^{-1*}$ restricts to define an $L$-linear isomorphism;\\

$Fr^{-1*}:L(C')\rightarrow L(C)^{1/p}$ $(\dag)$\\

We have that $Fr^{-1*}\circ Fr^{*}=Id$, restricted to $Rat_{k}$,
hence;\\

 $Fr^{-1*}\circ Fr^{*}:L(C)\rightarrow L(C')\rightarrow
L(C)^{1/p}$ $(\dag\dag)$\\

 is the inclusion map. Using the fact that $C_{p}$ and $C'$ are
 nonsingular projective curves, by $(\dag)$ we obtain an
 isomorphism $\theta:C_{p}\rightarrow C'$. By $(\dag\dag)$, we have
 that;\\

 $Fr\circ\theta=Frob:C_{p}\rightarrow C$\\

 Hence, without loss of generality, we can identify the morphisms
 $Fr$ and the more abstractly defined morphism $Frob$.

\end{rmk}

We now make the following further remark.

\begin{rmk}
Given the hypotheses of Theorem 2.1, with the modification that
$char(L)=p\neq 0$, we define a point $(ab)\in F$ to be wildly
ramified if $mult^{alg}_{(ab)}(F/D)$ is divisible by $p$. Theorem
2.1 holds excluding wildly ramified points, $(*)$. In order to see this,
we first replace the argument $(1)$, by showing that, for any
\emph{given} point $a\in D$, we can find a finite morphism $f$
from $D$ to $P^{1}(L)$, such that $f$ is etale in an open neighborhood
of $a$;\\

As $a$ is a non-singular, we can find a uniformising element $t$
in the local ring $O_{a,D}$ of $D$. Considering $t$ as an element
of the function field $L(D)$, we obtain an embedding $L(t)\subset
L(D)$, which, as $D$ is non-singular, determines a unique morphism
$f$ from $D$ to $P^{1}(L)$. Restricting the morphism to $A^{1}(L)$
and letting $x$ be the canonical coordinate, we have that
$f^{*}(x)=t$, hence $ord_{a}(f^{*}(x))=1$. This shows that $f$ is
etale in an open neighborhood of $a$ by Theorem 5.2 and Remarks
5.3 of \cite{depiro2}. $(\dag)$\\

As etale morphisms have multiplicity coprime to $p$, it is
sufficient to check the result $(*)$ for a branched cover
$\pi:F\rightarrow P^{1}(L)$. If $a\in F$ is not wildly ramified
for this cover, then we can follow through arguments $(2)$ and
$(3)$ of Theorem 2.1. The argument from $(4)$ to $(8)$ is the same
and we obtain the result of $(9)$ again using the fact that $m$
there is coprime to $p$. This proves the result $(*)$.\\

Theorem 2.1 also holds with the modification that $char(L)=p\neq
0$ and the cover $pr:F\rightarrow D$ is \emph{seperable}. However,
the proof requires more sophisticated methods, which we consider
below. We can, however, handle a special case by an
elementary counting argument. First observe that we can replace
the argument $(1)$ by observing that there exists a seperable
morphism $f$ from $D$ to $P^{1}(L)$. This either follows from the
argument $(\dag)$ above or using the classical result that the
function field $L(D)$ admits a seperating transcendence basis over
$L$, (see p27 of \cite{Hart}). Hence, it is sufficient to check
the result for a finite seperable cover $\pi:F\rightarrow
P^{1}(L)$. By a classical result, (see Proposition 2.2, p300, of
\cite{Hart}), there exist finitely many ramification points, in
particularly finitely many wild ramification points
$\{a_{1},\ldots,a_{n}\}$, for the cover $\pi$. By the previous
proof, we need only check the result of Theorem 2.1 for these
finitely many points. \\

Special Case. $a$ is a wild ramification point for the cover with
the property that that there exist no other wild ramification
points in the fibre $\pi^{-1}(\pi(a))$.\\

As both $F$ and $P^{1}(L)$ are non-singular, the finite morphism
$\pi$ is flat, by Lemma 5.11 of \cite{depiro2}. By a result in \cite{Mum},
(Corollary of Proposition 2, p218), we have that;\\

$\sum_{y\in\pi^{-1}(x)}mult_{y}^{alg}(F/P^{1})$ is independent of
$x\in P^{1}(L)$, and equals\\
\indent the cardinality of a generic fibre.\\

By Lemma 4.3 of \cite{depiro2}, a corresponding result also holds for Zariski
multiplicities. Hence, by the result of the previous proof in this
remark, the claim follows.\\

Unfortunately, one can have;\\

$a$ is a wild ramification point for the cover with the property
that there exist other wild ramification points
$\{a_{1},\ldots,a_{r}\}$, distinct from $a$, in the fibre
$\pi^{-1}(\pi(a))$.\\

It seems difficult to find any way of reducing this scenario to
the special case. However, one can still use a local method, which
is done in the following Theorem.

\end{rmk}

\begin{theorem}
Let hypotheses be as in Theorem 2.1, with the
modification that $char(L)=p\neq 0$ and the cover $pr:F\rightarrow
D$ is seperable. Then the notions of Zariski multiplicty and algebraic
multiplicity coincide.
\end{theorem}

\begin{proof}
By the previous remark, it is sufficient to consider the
case when $D$ is $P^{1}(L)$. Let $a\in F$, such that, without loss
of generality, $pr(a)=0$ in the restriction of $pr$ to $A^{1}(L)$.
As $a$ is non-singular, we can find polynomials
$\{f_{1},\ldots,f_{n-1}\}$ in the variables
$\{x_{1},\ldots,x_{n}\}$ of an affine coordinate system $A^{n}$,
such that $a$ corresponds to the origin $O$ of
this system and $F$ is defined locally by;\\

$f_{1}(x_{1},\ldots,x_{n})=\ldots=f_{n-1}(x_{1},\ldots,x_{n})=0$\\

with;\\

$Jac({f_{1},\ldots,f_{n-1}\over x_{2},\ldots,x_{n}})|_{\bar 0}\neq
0$\\

We may then apply the implicit function theorem, (see for example
p179 of \cite{Aby}), in order to find power series
$\{\eta_{1},\ldots,\eta_{n-1}\}$, in the variable $t$, with
$\eta_{j}(t)=0$, for $1\leq j\leq n-1$, such
that;\\

$f_{j}(t,\eta_{1}(t),\ldots,\eta_{n-1}(t))=0$, for $1\leq j\leq
n-1$. $(*)$\\

By $(*)$, we clearly have that the total transcendence degree of\\
$\{t,\eta_{1}(t),\ldots,\eta_{n-1}(t)\}$ over $L$ is equal to $1$.
Hence, we have that\\
 $\{\eta_{1}(t),\ldots,\eta_{n-1}(t)\}$ are
algebraic over $L(t)$. This implies, by the remarks at the
beginning of Section 3 of \cite{depiro3}, that they belong to the Henselisation of
$L[t]_{0}$, hence they define functions on some etale cover
$(U,0_{lift})$ with coordinate ring $L[t]^{ext}$ of $(A^{1},0)$. We have the ring map;\\

${L[x_{1},\ldots,x_{n}]\over
<f_{1},\ldots,f_{n-1}>}\rightarrow_{i} R=
{L[x_{1}]^{ext}[x_{2},\ldots,x_{n}]\over
<x_{2}-\eta(x_{1}),\ldots,x_{n}-\eta_{n-1}(x_{1})>}$\\

which corresponds to an etale cover $(U',a_{lift})$ of $(F,a)$. We
also have an isomorphism;\\

$R\rightarrow_{\gamma} L[t]^{ext}$; $x_{1}\mapsto t,x_{2}\mapsto
\eta_{1}(t),\ldots,x_{n}\mapsto
\eta_{n-1}(t)$\\

which corresponds to an isomorphism between $(U,0_{lift})$ and
$(U',a_{lift})$. Now consider the composition;\\

$\theta:(U,0_{lift})\rightarrow (U',a_{lift})\rightarrow
(F,a)\rightarrow_{pr}(A^{1},0)$\\

By the general method of \cite{depiro2}, we can define both the
algebraic and Zariski multiplicities of these covers. By
Theorem 1.4 and Lemma 2.2 of \cite{depiro3}, we have that;\\

$Mult_{(0,a)}(F/D)=Mult_{(0,0_{lift})}(U/A^{1})$\\

By Theorem 1.8 of \cite{depiro3}, we also have that;\\

$mult^{alg}_{(0,a)}(F/D)=mult^{alg}_{(0,0_{lift})}(U/A^{1})$\\

Hence, the theorem is shown by proving that Zariski multiplicity
and algebraic multiplicity coincide at $(0,0_{lift})$ for the
seperable cover $\theta$. Suppose that the algebraic multiplicity
is $m$, then, if $t$ is the canonical coordinate for $A^{1}$ at $0$, we have that;\\

$\theta^{*}t=t^{m}u(t)$\ for\ a\ unit\ $u(t)\in L[[t]]\cap
L(t)^{alg}$\\

By the usual factoring argument, see $(8)$ of Theorem 2.1,
it is sufficent to check that the Zariski multiplicity of
the seperable cover $\phi$ determined by;\\

$L[s]\rightarrow {L[t]^{ext}[s]\over <t^{m}u(t)-s>}$\\

is equal to $m$ at $(0,0_{lift})$ as well. This is done by the
general method of Lemmas 4.5 and 4.6 of \cite{depiro3}. We apply Weierstrass
preparation to $t^{m}u(t)-s$, see \cite{Aby} for the power
series version of this result, in order to obtain the factorisation;\\

$t^{m}u(t)-s=u(t,s)(t^{m}+c_{1}(s)t^{m-1}+\ldots+c_{m}(s))=u(t,s)g(t,s)$\\

where $c_{j}(s)\in L[[s]]\cap L(s)^{alg}$, $c_{j}(s)=0$ for $1\leq
j\leq m$ and $u(t,s)\in L[[s,t]]\cap L(s,t)^{alg}$ is a unit, see
Lemma 3.2 of \cite{depiro3}. As is done in Lemma 4.6 of \cite{depiro3}, we obtain
the etale cover determined by;\\

${L[t]^{ext}[s]\over <t^{m}u(t)-s>}\rightarrow {L[t,s]^{ext}\over
<u(t,s)g(t,s)>}$\\

By the argument there, it is sufficient to determine when the
Weierstrass factor $g(t,s)$ determines a generically reduced
cover. Using the method of resultants in Lemma 4.5 of \cite{depiro3},
this occurs if and only if ${\partial g\over\partial t}$ is not identically zero.
If ${\partial g\over\partial t}$ is identically zero, we obtain
the factorisation $g(t,s)=h(t^{p},s)$. This clearly implies that
the original cover $\phi$ is inseperable, which is a
contradiction. The theorem is then proved.

\end{proof}

We now have;\\

\begin{theorem}

Let hypotheses be as in Theorem 2.1, with the modification that
$char(L)=p\neq 0$. If $e$ denotes the Zariski multiplicity and $d$
the algebraic multiplicity at $a\in F$, then $d=ep^{n}$ and $\pi$
factors as $F\rightarrow_{h}F'\rightarrow_{g}D$ with $h=Frob^{n}$
and $g$ having algebraic multiplicity $e$ at $h(a)$.

\end{theorem}

\begin{proof}

As in Theorem 6.3 of \cite{depiro2}, we can factor $\pi$ into a purely inseperable
morphism $h:F\rightarrow F'$ and a seperable morphism
$g:F'\rightarrow D$  with $F'$ a smooth projective curve. We then
have a corresponding sequence of field extensions $L(D)\subset
L(F')\subset L(F)$, with $L(F)$ a purely inseperable extension of
$L(F')$. As $L(F)$ is a purely inseperable field extension of
$L(F')$, it has degree $p^{n}$ for some $n\geq 1$. Hence,
$L(F)=L(F')^{1/p^{n}}$ and we may, without loss of generality,
assume that $h=Frob^{n}$, see also Proposition 2.5 (p302) of
\cite{Hart}. By the previous theorem, the notions of Zariski multiplicity
and algebraic multiplicity coincide for the morphism $g$. By
Remarks 2.3, the Frobenius morphism $Frob$ may be identified with
$Fr$, without effecting Zariski or algebraic multiplicities.
Clearly, $Fr$ is a bijection on points, hence it is Zariski
unramified. $Fr$ has algebraic multiplicity $p$ everywhere, as,
for any point $x\in F'$, we can choose a local uniformiser
 $t$ at $x$ such that $Fr^{*}(t)=t^{p}$. It follows that $h$ has algebraic multiplicity
 $p^{n}$ everywhere and is Zariski unramified. The result now follows immediately from Lemma 4.5
  and Remarks 5.7 of \cite{depiro2}.\\

\end{proof}

We now give a local version of Theorem 2.1 in the general case of
algebraic curves over a field $L$ with $char(L)=0$ and find an
analogous version of Theorem 2.5, in the case when $char(L)=p\neq 0$.

\begin{theorem}
Let hypotheses be as in Theorem 3.3 of \cite{depiro2}, with the additional
assumption that $char(L)=0$ and $D$ is a \emph{smooth} curve.
Let $pr$ be the projection map of $F$ onto $D$. Then, if $(ab)\in F$ is
non-singular;\\

$Mult_{ab}(F/D)=mult_{ab}^{alg}(F/D)$\\

that is Zariski multiplicity and algebraic multiplicity coincide.
In particular, the cover $(F/D)$ is Zariski unramified at $(ab)$
iff there exists an open $U\subset F$, containing $(ab)$, such
that $pr:U\rightarrow D$ is etale.

\end{theorem}

\begin{proof}
For the first part of the theorem, we follow the proof of Theorem
2.1, the difference between the hypotheses there is that we do
\emph{not} assume that $F$ is smooth. Using the fact that
$D$ is smooth and the result of Theorem 2.1, we may, without loss
of generality, assume that $D=P^{1}(L)$. Now, one can follow
through the proof of Theorem 2.1, using the fact that $(ab)$ is
non-singular, in order to obtain the result. One should make the
modification that Zariski multiplicity is well defined for any
finite cover $F'\rightarrow F$ at $(abc)$ lying over $(ab)$. This
follows from an easy extension of Theorem 3.3 (in \cite{depiro2}),
to show that a nonsingular open subvariety of an irreducible projective variety
of dimension $1$ is presmooth (see \cite{Z}). For the second
part of the theorem, suppose that there exists an open $U\subset
F$, containing $(ab)$, such that $pr:U\rightarrow D$ is etale. As
$(ab)$ is non-singular, we may assume that $U$ defines a
non-singular open subvariety of $F$. Following the argument of
Theorem 2.1, from the end of $(4)$ to the end of $(7)$, we obtain
that the cover $(F/D)$ is Zariski unramified at $(ab)$. For the
converse, assume that the cover is Zariski unramified at $(ab)$.
By Theorem 5.2, Remarks 5.3 of \cite{depiro2} and the fact that $(ab)$ is
non-singular, it is sufficient to prove that
$d(pr):(m_{(ab)}/m_{(ab)}^{2})^{*}\rightarrow
(m_{a}/m_{a}^{2})^{*}$ is an isomorphism. Equivalently, we need to
show that the algebraic multiplicity $mult_{(ab)}^{alg}(F/D)$ of
$pr$ at $(ab)\in F$ equals $1$. This follows from the first part
of the theorem.

\end{proof}

\begin{theorem}
Let hypotheses be as in Theorem 3.3 of \cite{depiro2}, with the additional
assumption that $char(L)=p\neq 0$, $D$ is a \emph{smooth }curve
and the projection map $pr$ of $F$ onto $D$ is seperable. Then,
if $(ab)\in F$ is non-singular;\\

$Mult_{ab}(F/D)=mult_{ab}^{alg}(F/D)$\\

that is Zariski multiplicity and algebraic multiplicity coincide.
In particular, the cover $(F/D)$ is Zariski unramified at $(ab)$
iff there exists an open $U\subset F$, containing $(ab)$, such
that $pr:U\rightarrow D$ is etale.

\end{theorem}

\begin{proof}
Here, the hypotheses are the same as Theorem 2.5, with the modification that we do \emph{not} assume
$F$ is smooth. The proof is similar to the previous theorem. By Remarks 2.4, we
can assume that $D=P^{1}(L)$. Using the fact that $(ab)$ is
non-singular, one can either follow through the proof of Theorem
2.1, if $(ab)$ is not wildly ramified for the cover, or one can
use the method in Theorem 2.5, if $(ab)$ is
wildly ramified for the cover. For the second part, one can use the
 same reasoning as in the previous theorem.

\end{proof}

\begin{rmk}
This last result is required for the proof of Lemma 2.10 from \cite{depiro1} under suitable
assumptions, when $char(L)=p\neq 0$. The reader should consult the final section on Frobenius from
the paper \cite{depiro1}.
\end{rmk}

We finish this section with the following result;\\

\begin{theorem}
Let $G(X,Y)=0$ define an irreducible plane algebraic curve $C$,
with a non-singular point at $(0,0)$. Let $(T,\eta(T))$ be a power
series representation of this point. Then, for any plane, possibly
reduced, algebraic curve $F(X,Y)=0$ passing through $(0,0)$;\\

$F(T,\eta(T))\equiv 0$ iff $F$ contains $C$ as a component.\\

Otherwise, $I(G,F,(0,0))=ord_{T}F(T,\eta(T))$.

\end{theorem}

\begin{proof}
The proof partly uses the methods of \cite{depiro3}. For
the first part, note that if $F$ contains $C$ as a component, then
by the Nullstellenstatz, there exists $H(X,Y)$ such that
$F(X,Y)=H(X,Y)G(X,Y)$. It then follows trivially that
$F(T,\eta(T))\equiv 0$. For the converse direction, suppose that
$F(T,\eta(T))\equiv 0$. As in Lemma 4.17 of \cite{depiro3}, we may interpret the
equation $Y-\eta(X)$ as defining a curve $C_{1}$ on some etale
extension $i:(A^{2}_{et},(00)^{lift})\rightarrow (A^{2},(00))$
such that $i(C_{1})\subset C$. The vanishing of $F(X,Y)$ on
$C_{1}$ then implies that $F$ intersects $C$ in an open dense
subset. Therefore, as both $F$ and $C$ define Zariski closed sets,
$F$ must contain $C$ as a component. For the second part of the
theorem, we may therefore assume that $F$ has finite intersection
with $C$ and $ord_{T}F(T,\eta(T))$ is defined. Suppose that
$F(X,Y)$ has degree $d$ and consider $F$ as part of the family of
degree $d$ curves $Q_{d}$. Without loss of generality, we may
suppose that $F(X,Y)=H(X,Y,\bar v^{0})$ where, for $\bar v\in
Par_{Q_{d}}$, $H(X,Y,\bar v)$ defines an algebraic curve of degree
$d$. Similarily, we can write $G(X,Y)$ in the form $G(X,Y,\bar
u^{0})$ for some non-varying constant $\bar u^{0}$.
As in Lemma 4.17 of \cite{depiro3}, we have the sequence of maps;\\

$L[\bar v]\rightarrow {L[X,Y][\bar v]\over <G(X,Y,\bar
u_{0}),H(X,Y,\bar v)>}\rightarrow {L[X]^{ext}[Y][\bar v]\over
<Y-\eta(X),H(X,Y,\bar
v)>}$\\

which corresponds to a sequence of finite covers;\\

$F_{1}\rightarrow F'(\bar u^{0},V)\rightarrow Spec(L[\bar v])$\\

One checks that the left hand morphism is etale at $(\bar
v^{0},(00)^{lift})$, by direct calculation. We use the fact that
$F$ is non-singular at $(00)$, therefore the completion of the
local rings ${L[X,Y]\over <G(X,Y,\bar u_{0})>}_{(00)}$ and
${L[X]^{ext}[Y]\over <Y-\eta(X)>}_{(00)}$ are in both cases equal
to the formal power series ring $L[[X]]$.\\

We now compute the Zariski multiplicity of the cover
$F_{1}\rightarrow Spec(L[\bar v])$ at $(\bar v^{0},(00)^{lift})$
$(*)$. We are given the formal power series $H(X,\eta(X),\bar
v)\in L[[X,\bar v]]$. Let $d=ord_{X}H(X,\eta(X),\bar v_{0})$.
Then, by Weierstrass preparation in several variables, see
\cite{Aby}, we can find $H_{1}(X,\bar v)$
and $U(X,\bar v)$ in $L[[X,\bar v]]$ such that;\\

$H(X,\eta(X),\bar v)=H_{1}(X,\bar v)U(X,\bar v)$\\

and $U(0,\bar v_{0})\neq 0$ and\\

$H_{1}(X,\bar v)=X^{d}+c_{1}(\bar v)X^{d-1}+\ldots c_{d}(\bar
v)$\\

with $c_{j}(\bar v_{0})=0$ for $1\leq j\leq d$. Now use the proofs
of Lemma 4.5 and 4.6 from \cite{depiro3} and the fact that the cover;\\

$Spec(H_{1}(X,\bar v))\rightarrow Spec(L[\bar v])$\\

is generically reduced to show the Zariski multiplicity of the
cover $(*)$ is exactly $d$. This proves that the Zariski
multiplicity of the cover;\\

$F'(\bar u^{0},V)\rightarrow Spec(L[\bar v])$\\

at $((0,0),\bar v^{0})$ is exactly $d$ as well. By the general
result of the paper \cite{depiro3}, that;\\

$I(C_{\bar u^{0}},C_{\bar v^{0}},(00))=RightMult_{(00)}(C_{\bar
u^{0}},C_{\bar v^{0}})$\\

when $C_{\bar u^{0}}$ defines a reduced algebraic curve, the
result of the theorem follows.

\end{proof}

\begin{rmk}
This last Theorem was required in the proof of Theorem 6.1 of \cite{depiro1}. It is also required in the proof of Remarks 4.8 below.
\end{rmk}

\end{section}
\begin{section}{A refined theory of $g_{n}^{r}$}

The purpose of this section is to refine the general theory of
$g_{n}^{r}$, given in \cite{depiro1}, in order to take into account the notion of a branch
for a projective algebraic curve. We will rely heavily on results
proved in \cite{depiro1}. We also refer the reader there for the
relevant notation. We will make \emph{no} assumptions on
the characteristic of the base field $L$. As usual, by an algebraic curve, we always mean
a projective irreducible variety of dimension $1$.\\

\begin{defn}

Let $C\subset P^{w}$ be a projective algebraic curve of degree $d$
and let $\Sigma$ be a linear system of dimension $R$, contained in
the space of algebraic forms of degree $e$ on $P^{w}$. Let
$\phi_{\lambda}$ belong to $\Sigma$, having finite intersection
with $C$. Then, if $p\in C\cap\phi_{\lambda}$ and $\gamma_{p}$ is
a branch centred at $p$, we define;\\

$I_{p}(C,\phi_{\lambda})=I_{italian}(p,C,\phi_{\lambda})$\\

$I_{p}^{\Sigma}(C,\phi_{\lambda})=I_{italian}^{\Sigma}(p,C,\phi_{\lambda})$\\

$I_{p}^{\Sigma,mobile}(C,\phi_{\lambda})=I_{italian}^{\Sigma,mobile}(p,C,\phi_{\lambda})$\\

$I_{\gamma_{p}}(C,\phi_{\lambda})=I_{italian}(p,\gamma_{p},C,\phi_{\lambda})$\\

$I_{\gamma_{p}}^{\Sigma}(C,\phi_{\lambda})=I_{italian}^{\Sigma}(p,\gamma_{p},C,\phi_{\lambda})$\\

$I_{\gamma_{p}}^{\Sigma,mobile}(C,\phi_{\lambda})=I_{italian}^{\Sigma,mobile}(p,\gamma_{p},C,\phi_{\lambda})$\\

where $I_{italian}$ was defined in \cite{depiro1}.

\end{defn}

It follows that, as $\lambda$ varies in $Par_{\Sigma}$, we obtain
a series of weighted sets;\\

$W_{\lambda}=\{n_{\gamma_{p_{1}}^{1}},\ldots,n_{\gamma_{p_{1}}^{n_{1}}},\ldots,n_{\gamma_{p_{m}}^{1}},\ldots,n_{\gamma_{p_{m}}^{n_{m}}}\}$\\

where;\\

$\{p_{1},\ldots,p_{i},\ldots,p_{m}\}=C\cap\phi_{\lambda}$, for $1\leq i\leq m$,\\

$\{\gamma_{p_{i}}^{1},\ldots,\gamma_{p_{i}}^{j(i)},\ldots,\gamma_{p_{i}}^{n_{i}}\}$,
for $1\leq j(i)\leq n_{i}$, consists of the branches of $C$ centred at $p_{i}$\\

and\\

$I_{\gamma_{p_{i}}^{j(i)}}(C,\phi_{\lambda})=n_{\gamma_{p_{i}}^{j(i)}}$\\

By the branched version of the Hyperspatial Bezout Theorem, see
\cite{depiro1}, the total weight of any of these sets, which we
will occasionally abbreviate by $C\sqcap\phi_{\lambda}$, is always
equal to $de$. Let $r$ be the least integer such that every
weighted set $W_{\lambda}$ is defined by a linear subsystem
$\Sigma'\subset\Sigma$ of dimension $r$.\\

\begin{defn}

We define;\\

$Series(\Sigma)=\{W_{\lambda}:\lambda\in Par_{\Sigma}\}$\\

$dimension(Series(\Sigma))=r$\\

$order(Series(\Sigma))=de$\\

\end{defn}

We then claim the following;\\

\begin{theorem}.\\

(i). $r\leq R$, with equality iff every weighted set $W_{\lambda}$
of the series is cut out by a \emph{single} form of $\Sigma$.\\

(ii). $r\lneqq R$ iff there exists a form $\phi_{\lambda}$ in
$\Sigma$, containing all of $C$.
\end{theorem}

\begin{proof}
We first show the equivalence of $(i)$ and $(ii)$. Suppose that
$(i)$ holds and $r\lneqq R$. Then, we can find a weighted set $W$
and distinct elements $\{\lambda_{1},\lambda_{2}\}$ of
$Par_{\Sigma}$ such that $W=W_{\lambda_{1}}=W_{\lambda_{2}}$. Let
$\{\phi_{\lambda_{1}},\phi_{\lambda_{2}}\}$ be the corresponding
algebraic forms of $\Sigma$ and consider the pencil
$\Sigma_{1}\subset \Sigma$ defined by these forms. We claim
that;\\

$W=C\sqcap(\mu_{1}\phi_{\lambda_{1}}+\mu_{2}\phi_{\lambda_{2}})$,
for $[\mu_{1}:\mu_{2}]\in P^{1}$ $(*)$\\

This follows immediately from the results in \cite{depiro1} that
the condition of multiplicity at a branch is \emph{linear} and the
branched version of the Hyperspatial Bezout Theorem. Now choose a
point $p\in C$, which is not a base point for any of the branches in $W$.
Then, the condition that an
algebraic form $\phi_{\lambda}$ passes through $p$ defines a
hyperplane condition on $Par_{e}$, hence, intersects
$Par_{\Sigma_{1}}$ in a point. Let $\phi_{\lambda_{0}}$ be the
algebraic form in $\Sigma_{1}$ defined by this parameter. Then, by
$(*)$, we have that;\\

$W\cup\{p\}\subseteq C\sqcap\phi_{\lambda_{0}}$\\

Hence, the total multiplicity of intersection of
$\phi_{\lambda_{0}}$ with $C$ is at least equal to $de+1$. By the
branched version of the Hyperspatial Bezout Theorem, $C$ must be
contained in $\phi_{\lambda_{0}}$. Conversely, suppose that $(i)$
holds and there exists a form $\phi_{\lambda_{0}}$ in $\Sigma$
containing all of $C$. Let $W$ be cut out by $\phi_{\lambda_{1}}$
and consider the pencil $\Sigma_{1}\subset\Sigma$ generated by
$\{\phi_{\lambda_{0}},\phi_{\lambda_{1}}\}$. By the same argument
as above, we can find $\phi_{\lambda_{2}}$ in $\Sigma_{1}$,
distinct from $\phi_{\lambda_{1}}$, which also cuts out $W$.
Hence, by $(i)$, we must have that $r\lneq R$. Therefore, $(ii)$
holds. \\

The argument that $(ii)$ implies $(i)$ is similar.\\

We now prove that $(ii)$ holds. Using the Hyperspatial Bezout
Theorem, the condition on $Par_{\Sigma}$ that a form
$\phi_{\lambda}$ contains $C$ is linear. Let $H$ be the linear
subsystem of $\Sigma$, consisting of forms containing $C$ and let
$h=dim(H)$. Let $K\subset\Sigma$ be a maximal linear subsystem,
having finite intersection with $C$. Then $K$ has no form in
common with $H$ and $dim(K)=R-h-1$. We claim that every weighted
set in $Series(\Sigma)$ is cut out by a unique form from $K$. For
suppose that $W=C\sqcap\phi_{\lambda}$ is such a weighted set and
consider the linear system defined by $<H,\phi_{\lambda}>$. If
$\phi_{\mu}$ belongs to this system and has finite intersection
with $C$, then clearly $(C\cap\phi_{\lambda})=(C\cap\phi_{\mu})$.
Using linearity of multiplicity at a branch and the Hyperspatial
Bezout Theorem again (by convention, a form containing $C$ has
infinite multiplicity at a branch), we must have that
$(C\sqcap\phi_{\lambda})=(C\sqcap\phi_{\mu})$. Now consider
$K\cap <H,\phi_{\lambda}>$. We have that;\\

$codim(K\cap <H,\phi_{\lambda}>)\leq
codim(K)+codim(<H,\phi_{\lambda}>)$\\
$. \ \ \ \ \ \ \ \ \ \ \ \ \ \ \ \ \ \ \ \ \ \ \ \ \ \ \ \ \ \ \ \ =(h+1)+(R-(h+1))=R$.\\

 Hence,
$dim(K\cap <H,\phi_{\lambda}>)\geq 0$. We can, therefore, find a
form $\phi_{\mu}$ belonging to $K$ such that
$W=(C\sqcap\phi_{\mu})$. We need to show that $\phi_{\mu}$ is the
unique form in $K$ defining $W$. This follows by the argument
given above. It follows immediately that $r=dim(K)=R-h-1$. Hence,
$r\lneq R$ iff $h\geq 0$. Therefore, $(ii)$ is shown.

\end{proof}

Using this theorem, we give a more refined definition of a
$g_{n}^{r}$.\\

\begin{defn}
Let $C\subset P^{w}$ be a projective algebraic curve. By a
$g_{n}^{r}$ on $C$, we mean the collection of weighted sets,
without repetitions, defined by $Series(\Sigma)$ for \emph{some}
linear system $\Sigma$, such that $r=dimension(Series(\Sigma))$
and $n=order(Series(\Sigma))$. If a branch $\gamma_{p}^{j}$
appears with multiplicity at least $s$ in every weighted set of a
$g_{n}^{r}$, as just defined, then we allow the possibility of
removing some multiplicity contribution $s'\leq s$ from each
weighted set and adjusting $n$ to $n'=n-s'$.

\end{defn}

\begin{rmk}
The reader should observe carefully that a $g_{n}^{r}$ is defined
independently of a particular linear system. However, by the
previous theorem, for any $g_{n}^{r}$, there exists a $g_{n'}^{r}$
with $n\leq n'$ such that the following property holds. The
$g_{n'}^{r}$ is defined by a linear system of dimension $r$,
having finite intersection with $C$, such that each there is a
bijection between the weighted sets $W$ in the $g_{n'}^{r}$ and
the $W_{\lambda}$ in $Series(\Sigma)$. The original $g_{n}^{r}$ is
obtained from the $g_{n'}^{r}$ by removing some fixed branch
contribution.
\end{rmk}

We now reformulate the results of Section 2 and Section 5 in
\cite{depiro1} for this new definition of a $g_{n}^{r}$. In order
to do this, we require the following definition;\\

\begin{defn}
Suppose that $C\subset P^{w}(L)$ is a projective algebraic curve
and $C^{ext}\subset P^{w}(K)$ is its non-standard model. Let a
$g_{n}^{r}$ be given on $C$, defined by a linear system $\Sigma$
after removing some fixed branch contribution. We define the
extension $g_{n}^{r,ext}$ of the $g_{n}^{r}$ to the nonstandard
model $C^{ext}$ to be the collection of weighted sets, without
repetitions, defined by $Series(\Sigma)$ on $C^{ext}$, after
removing the same fixed point contribution. Note that, by
definability of multiplicity at a branch, see Theorem 6.5 of
\cite{depiro1}, if $\gamma_{p}^{j}$ is a branch of $C$ and;\\

$I_{italian}(p,\gamma_{p}^{j},C,\phi_{\lambda})\geq k$,
$(\lambda\in Par_{\Sigma(L)})$\\

then;\\

$I_{italian}(p,\gamma_{p}^{j},C,\phi_{\lambda})\geq k$,
$(\lambda\in Par_{\Sigma(K)})$\\

Hence, it \emph{is} possible to remove the same fixed point
contribution
of $Series(\Sigma)$ on $C^{ext}$. See also the proof of Lemma 3.7.\\

\end{defn}

It is a remarkable fact that, after introducing the notion of a
branch, the definition is independent of the particular linear
system $\Sigma$. This is the content of the following lemma;\\

\begin{lemma}
The previous definition is independent of the particular choice of
linear system $\Sigma$ defining the $g_{n}^{r}$.
\end{lemma}

\begin{proof}

We divide the proof into the following cases;\\

Case 1. $\Sigma\subset\Sigma'$;\\

By the proof of Theorem 3.3, we can find a linear system
$\Sigma_{0}\subset\Sigma\subset\Sigma'$ of dimension $r$, having
finite intersection with $C$, such that the $g_{n}^{r}$ is defined
by removing some fixed contribution from $\Sigma_{0}$. Here, we
have also used the fact that the base point contributions (at a
branch) of $\{\Sigma_{0},\Sigma,\Sigma'\}$ are the same. Again, by
Theorem 3.3, if $W_{\lambda'}$ is a weighted set defined by
$\Sigma'$ on $C^{ext}$, then it appears as a weighted set
$V_{\lambda''}$ defined by $\Sigma_{0}$ on $C^{ext}$. Hence, it
appears as a weighted set $V_{\lambda''}$ defined by $\Sigma$ on
$C^{ext}$. By the converse argument and Remarks 3.5 on base branch
contributions, the proof is shown.\\

Case 2. $\Sigma$ are $\Sigma'$ are both linear systems of
dimension $r$, having finite intersection with $C$, such that $degree(\Sigma)=degree(\Sigma')=n$;\\

By Theorem 3.3, every weighted set $W$ in the $g_{n}^{r}$ is
defined uniquely by weighted sets $W_{\lambda_{1}}$ and
$V_{\lambda_{2}}$ in $Series(\Sigma_{1})$ and $Series(\Sigma_{2})$
respectively. Let $(C^{ns},\Phi^{ns})$ be a non-singular model of
$C$. Using the method of Section 5 in \cite{depiro1} to avoid the
technical problem of presentations of $\Phi^{ns}$ and base point
contributions, we may, without loss of generality, assume that
there exist finite covers $W_{1}\subset Par_{\Sigma}\times C^{ns}$ and $W_{2}\subset Par_{\Sigma'}\times C^{ns}$ such that;\\

$j_{k,\Sigma}(\lambda,p_{j})\equiv
Mult_{(W_{1}/Par_{\Sigma})}(\lambda,p_{j})\geq k$ iff
$I_{italian}(p,\gamma_{p}^{j},C,\phi_{\lambda})\geq k$\\

$j_{k,\Sigma'}(\lambda',p_{j})\equiv
Mult_{(W_{2}/Par_{\Sigma'})}(\lambda',p_{j})\geq k$ iff
$I_{italian}(p,\gamma_{p}^{j},C,\psi_{\lambda'})\geq k$\\

Then consider the sentences;\\

$(\forall \lambda\in Par_{\Sigma})(\exists!\lambda'\in
Par_{\Sigma'})\forall x\in
C^{ns}[\bigwedge_{k=1}^{n}(j_{k}(\lambda,x)\leftrightarrow
j_{k}(\lambda',x))]$\\

$(\forall \lambda'\in Par_{\Sigma})(\exists!\lambda\in
Par_{\Sigma})\forall x\in
C^{ns}[\bigwedge_{k=1}^{n}(j_{k}(\lambda',x)\leftrightarrow
j_{k}(\lambda,x))]$ (*)\\

in the language of $<P^{1}(L),C_{i}>$, considered as a Zariski
structure with predicates $\{C_{i}\}$ for Zariski closed subsets
defined over $L$, (see \cite{Z}). We have, again by results of
\cite{Z} or \cite{depiro4}, that $<P^{1}(L),C_{i}>\prec
<P^{1}(K),C_{i}>$, for the nonstandard model $P(K)$ of $P(L)$. It
follows immediately from the algebraic definition of $j_{k}$ in
\cite{Z}, that, for any weighted set $W_{\lambda_{1}'}$ defined by
$Series(\Sigma)$ on $C^{ext}$, there exists a unique weighted set
$V_{\lambda_{2}'}$ defined by $Series(\Sigma')$ on $C^{ext}$ such
that $W_{\lambda_{1}'}=V_{\lambda_{2}'}$, and conversely. Hence,
the proof is shown.\\

Case 3. $\Sigma$ are $\Sigma'$ are both linear systems of
dimension $r$, having finite intersection with $C$;\\

Let $n_{1}=degree(\Sigma)$ and $n_{2}=degree(\Sigma')$. Then the
original $g_{n}^{r}$ is obtained from $Series(\Sigma)$, by
removing a fixed point contribution of multiplicity $n_{1}-n$,
and, is obtained from $Series(\Sigma')$, by removing a fixed point
contribution of multiplicity $n_{2}-n$. We now imitate the proof
of Case 2, with the slight modification that, in the construction
of the sentences given by $(*)$, we make an adjustment of the
multiplicity statement at the finite number of branches where a
fixed point contribution has been removed. The details are left to
the reader.
\end{proof}

Now, using Definition 3.6, we construct a specialisation operator
$sp:g_{n}^{r,ext}\rightarrow g_{n}^{r}$. We first require the following simple
lemma;\\

\begin{lemma}
Let $C\subset P^{w}(L)$ be a projective algebraic curve and let
$C^{ext}\subset P^{w}(K)$ be its nonstandard model. Let $p'\in
C^{ext}$ be a non-singular point, with specialisation $p\in C$.
Then there exists a unique branch $\gamma_{p}^{j}$ such that
$p'\in\gamma_{p}^{j}$.

\end{lemma}

\begin{proof}
We may assume that $p'\neq p$, otherwise $p$ would be non-singular
and, by Lemma 5.4 of \cite{depiro1}, would be the origin of a
single branch $\gamma_{p}$. Let $(C^{ns},\Phi)$ be a non-singular
model of $C$, then $p'$ must belong to the canonical set
$V_{[\Phi]}$, hence there exists a unique $p''\in C^{ns}$ such
that $\Phi(p'')=p'$. By properties of specialisations, $p''\in
C^{ns}\cap{\mathcal V}_{p_{j}}$ for some
$p_{j}\in\Gamma_{[\Phi]}(x,p)$. Hence, by definition of a branch
given in Definition 5.15 of \cite{depiro1}, we must have that
$p'\in\gamma_{p}^{j}$. The uniqueness statement follows as well.
\end{proof}

We now make the following definition;\\

\begin{defn}
Let $C\subset P^{w}(L)$ be a projective algebraic curve and let
$C^{ext}\subset P^{w}(K)$ be its non-standard model. Given a
$g_{n}^{r}$ on $C$ with extension $g_{n}^{r,ext}$ on $C^{ext}$,
we define the specialisation operator;\\

$sp:g_{n}^{r,ext}\rightarrow g_{n}^{r}$\\

by;\\

$sp(\gamma_{p'})=\gamma_{p}^{j}$, for $p'\in NonSing(C^{ext})$ and $\gamma_{p}^{j}$ as in Lemma 3.8.\\

$sp(\gamma_{p}^{j})=\gamma_{p}^{j}$, for $p\in
Sing(C^{ext})=Sing(C)$ and
$\{\gamma_{p}^{1},\ldots,\gamma_{p}^{j},\ldots,\gamma_{p}^{s}\}$
\indent\indent\indent\indent\indent\indent\  enumerating the branches at $p$.\\

$sp(n_{1}\gamma_{p_{1}}^{j_{1}}+\ldots+n_{r}\gamma_{p_{r}}^{j_{r}})=n_{1}sp(\gamma_{p_{1}}^{j_{1}})+\ldots+n_{r}sp(\gamma_{p_{r}}^{j_{r}})$,\\

for a linear combination of branches with $n_{1}+\ldots+n_{r}=n$\\

\end{defn}

It is also a remarkable fact that, after introducing the notion of
a branch, the specialisation operator $sp$ is well defined. This is the content of the following lemma;\\

\begin{lemma}
Let hypotheses be as in the previous definition, then, if $W$ is a
weighted set belonging to $g_{n}^{r,ext}$, its specialisation
$sp(W)$ belongs to $g_{n}^{r}$.

\end{lemma}

\begin{proof}
We may assume that there exists a linear system $\Sigma$, having
finite intersection with $C$, such that $dimension(\Sigma)=r$ and
$degree(\Sigma)=n_{1}$, with the $g_{n}^{r}$ and $g_{n}^{r,ext}$
both defined by $Series(\Sigma)$, after removing some fixed branch
contribution $W_{0}$ of multiplicity $n_{1}-n$. Let $W$ be a
weighted set of the $g_{n}^{r,ext}$, then $W\cup
W_{0}=(C\sqcap\phi_{\lambda'})$, for some unique $\lambda'\in
Par_{\Sigma}$. We claim that $sp(W\cup
W_{0})=C\sqcap\phi_{\lambda}$, for the specialisation $\lambda\in
Par_{\Sigma}$ of $\lambda'$ $(*)$. As $sp(W_{0})=W_{0}$, it then
follows immediately from linearity of $sp$, that $sp(W)$ belongs
to the $g_{n}^{r}$ as required. We now show $(*)$. Let $p\in C$ and
let $\gamma_{p}$ be a branch centred at $p$. By $\gamma_{p}^{ext}$, we mean the branch
at $p$, where $p$ is considered as an element of $C^{ext}$. We now claim that;\\

$I_{\gamma_{p}}(C,\phi_{\lambda})=I_{\gamma_{p}^{ext}}(C,\phi_{\lambda'})+\sum_{p'\in(\gamma_{p}\setminus
p)}I_{\gamma_{p'}^{ext}}(C,\phi_{\lambda'})$ $(**)$\\

Let $(C^{ns},\Phi)\subset P^{w'}(L)$ be a non-singular model of
$C$, such that $\gamma_{p}$ corresponds to $C^{ns,ext}\cap
{\mathcal V}_{q}$, where $q\in\Gamma_{[\Phi]}(x,p)$ and ${\mathcal
V}_{q}$ is defined relative to the specialisation from $P(K)$ to
$P(L)$. Let $C^{ns,ext,ext}\subset P^{w'}(K')$ be a non-standard
model of $C^{ns,ext}$, such that $\gamma_{q}^{ext}$ corresponds to
$C^{ns,ext,ext}\cap {\mathcal V}_{q}$, where ${\mathcal V}_{q}$ is
defined relative to the specialisation from $P(K')$ to $P(K)$.
Then, for $p'\in {(\gamma_{p}\setminus p)}$, we can find $q'\in
{\mathcal V}_{q}\cap C^{ns,ext}$ such that $\gamma_{p'}$
corresponds to ${\mathcal V}_{q'}\cap C^{ns,ext,ext}$. We may
choose a suitable presentation $\Phi_{\Sigma_{1}}$ of $\Phi$, such
that $Base(\Sigma_{1})$ is disjoint from $\Gamma_{[\Phi]}(x,p)$,
and, therefore, disjoint from $\Gamma_{[\Phi]}(x,p')$, for
$p'\in{(\gamma_{p}\setminus p)}$. Let
$\{\overline{\phi_{\lambda}}\}$ denote the lifted family of on
$C^{ns}$ from the presentation $\Phi_{\Sigma'}$. In this case,
we have, by results of \cite{depiro1}, that;\\

$I_{\gamma_{p}}(C,\phi_{\lambda})=I_{q}(C^{ns},\overline{\phi_{\lambda}})$\\

$I_{\gamma_{p}^{ext}}(C,\phi_{\lambda'})=I_{q}(C^{ns},\overline{\phi_{\lambda'}})$\\

$I_{\gamma_{p'}^{ext}}(C,\phi_{\lambda'})=I_{q'}(C^{ns},\overline{\phi_{\lambda'}})$ $(1)$\\

By summability of specialisation, see \cite{depiro1} and \cite{depiro3};\\

$I_{q}(C^{ns},\overline{\phi_{\lambda}})=I_{q}(C^{ns},\overline{\phi_{\lambda'}})+\sum_{q'\in
C^{ns}\cap{({\mathcal
V}_{q}}\setminus q)}I_{q'}(C^{ns},\overline{\phi_{\lambda'}})$ $(2)$\\

Combining $(1)$ and $(2)$, the result $(**)$ follows, as required.
Now, suppose that a branch $\gamma_{p}$ occurs with non-trivial
multiplicity in\\ $sp(C\sqcap\phi_{\lambda'})$. By Definition 3.9,
the contribution must come from either
$I_{\gamma_{p}^{ext}}(C,\phi_{\lambda'})$ or
$I_{\gamma_{p'}^{ext}}(C,\phi_{\lambda'})$, for some
$p'\in({\gamma_{p}\setminus p})$. Applying $sp$ to $(**)$, one
sees that the branch $\gamma_{p}$ occurs with multiplicity
$I_{\gamma_{p}}(C,\phi_{\lambda})$. It follows that
$sp(C\sqcap\phi_{\lambda'})=C\sqcap\phi_{\lambda}$, hence $(*)$ is
shown. The lemma then follows.

\end{proof}

We can now reformulate the results of Section 2 and Section 5 of
\cite{depiro1} in the language of this refined theory of $g_{n}^{r}$. We first make the
following definition;\\

\begin{defn}
Let $C\subset P^{w}$ be a projective algebraic curve and let a
$g_{n}^{r}$ be given on $C$. Let $W$ be a weighted set in this
$g_{n}^{r}$ or its extension $g_{n}^{r,ext}$ and let $\gamma_{p}$
be a branch centred at $p$. Then we say
that;\\

$\gamma_{p}$ is $s$-fold ($s$-plo) for $W$ if it appears with
multiplicity at least $s$.\\

$\gamma_{p}$ is multiple for $W$ if it appears with multiplicity
at least $2$.\\

$\gamma_{p}$ is simple for $W$ if it is not multiple.\\

$\gamma_{p}$ is counted (contato) $s$-times in $W$ if it appears
with multiplicity exactly $s$.\\

$\gamma_{p}$ is a base branch of the $g_{n}^{r}$ if it
appears in \emph{every} weighted set.\\

$\gamma_{p}$ is $s$-fold for the $g_{n}^{r}$ if it is $s$-fold in
$W$ for \emph{every} weighted set $W$ of the $g_{n}^{r}$.\\

$\gamma_{p}$ is counted $s$-times for the $g_{n}^{r}$ if it is
$s$-fold for the $g_{n}^{r}$ and is counted $s$-times in some
weighted set $W$ of the $g_{n}^{r}$.
\end{defn}

We then have the following;\\

\begin{theorem}{Local Behaviour of a $g_{n}^{r}$}\\

Let $C$ be a projective algebraic curve and let a $g_{n}^{r}$ be
given on $C$. Let $\gamma_{p}$ be a branch centred at $p$, such
that $\gamma_{p}$ is counted $s$-times for the $g_{n}^{r}$. If
$\gamma_{p}$ is counted $t$ times in a given weighted set $W$,
then there exists a weighted set $W'$ in $g_{n}^{r,ext}$ such that
$sp(W')=W$ and $sp^{-1}(t\gamma_{p})$ consists of the branch
$\gamma_{p}$ counted $s$-times and $t-s$ other distinct branches
$\{\gamma_{p_{1}},\ldots,\gamma_{p_{t-s}}\}$, each counted once in
$W'$.
\end{theorem}

\begin{proof}
Without loss of generality, we may assume that the $g_{n}^{r}$ is
defined by a linear system $\Sigma$ of dimension $r$, having
finite intersection with $C$. Let $W$ be the weighted set defined
by $\phi_{\lambda}$ in $\Sigma$. Suppose that $s=0$, then
$\gamma_{p}$ is not a base branch for $\Sigma$. Hence, by Lemma
5.25 of \cite{depiro1}, we can find $\lambda'\in{\mathcal
V}_{\lambda}$, generic in $Par_{\Sigma}$, and distinct
$\{p_{1},\ldots,p_{t}\}=C^{ext}\cap\phi_{\lambda'}\cap
({\gamma_{p}\setminus p})$ such that the intersections at these
points are transverse. Let $W'$ be the weighted set defined by
$\phi_{\lambda'}$ in $g_{n}^{r,ext}$. By the proof of $(*)$ in
Lemma 3.10, we have that $sp(W')=W$. By the construction of $sp$
in Definition 3.9, we have that $sp^{-1}(t\gamma_{p})$ consists of
the distinct branches $\{\gamma_{p_{1}},\ldots,\gamma_{p_{t}}\}$,
each counted once in $W'$. If $s\geq 1$, then $\gamma_{p}$ is a
base branch for $\Sigma$. By Lemma 5.27 of \cite{depiro1}, we have
that
$I^{\Sigma,mobile}_{italian}(p,\gamma_{p},C,\phi_{\lambda})=t-s$.
The result then follows by application of Lemma 5.28 in
\cite{depiro1} and the argument given above.
\end{proof}

We now note the following;\\

\begin{lemma}
Let a $g_{n}^{r}$ be given on a projective algebraic curve $C$.
Let $W_{0}$ be \emph{any} weighted set on $C$ with total
multiplicity $n'$. Then the collection of weighted sets given by
$\{W\cup W_{0}\}$ for the weighted sets $W$ in the $g_{n}^{r}$
defines a $g_{n+n'}^{r}$.

\end{lemma}

\begin{proof}
Let the original $g_{n}^{r}$ be obtained from a linear system
$\Sigma$ of dimension $r$ and degree $n''$, having finite
intersection with $C$, after removing some fixed branch
contribution $J$ of total multiplicity $n''-n$. Let
$\{\phi_{0},\ldots,\phi_{r}\}$ be a basis for $\Sigma$ and let
$\{n_{1}\gamma_{p_{1}}^{j_{1}},\ldots,n_{m}\gamma_{p_{m}}^{j_{m}}\}$
be the branches appearing in $W_{0}$ with total multiplicity
$n_{1}+\ldots+n_{m}=n'$ $(\dag)$. Let $\{H_{1},\ldots,H_{m}\}$ be
hyperplanes passing through the points $\{p_{1},\ldots,p_{m}\}$
and let $G$ be the algebraic form of degree $n'$ defined by
$H_{1}^{n_{1}}\centerdot\ldots\centerdot H_{m}^{n_{m}}$. Let
$\Sigma'$ be the linear system of dimension $r$ defined by the
basis\\
 $\{G\centerdot\phi_{0},\ldots,G\centerdot\phi_{r}\}$. As
we may assume that $C$ is not contained in any hyperplane section,
$\Sigma'$ has finite intersection with $C$. We claim that
$g_{n''}^{r}(\Sigma)\subset g_{n''+n'deg(C)}^{r}(\Sigma')$, in the
sense that every weighted set $W_{\lambda}$ defined by
$g_{n''+n'deg(C)}^{r}(\Sigma')$ is obtained from the corresponding
$V_{\lambda}$ in $g_{n''}^{r}(\Sigma)$ by adding a \emph{fixed}
weighted set $W_{1}\supset W_{0}$ of total multiplicity $n'deg(C)$
$(*)$. The proof then follows as we can recover the original
$g_{n}^{r}$ by removing the fixed branch contribution
$J\cup(W_{1}\setminus W_{0})$ from
$g_{n''+n'deg(C)}^{r}(\Sigma')$. In order to show $(*)$, let
$W_{1}$ be the weighted set defined by $C\sqcap G$. By the
branched version of the Hyperspatial Bezout Theorem, see Theorem
5.13 of \cite{depiro1}, this has total multiplicity $n'deg(C)$. We
claim that $W_{0}\subset W_{1}$ $(**)$. Let $\gamma_{p}^{j}$ be a
branch appearing in $(\dag)$ with multiplicity $s$. By
construction, we can factor $G$ as $H^{s}\centerdot R$, where $H$
is a hyperplane
passing through $s$. We need to show that;\\

$I_{\gamma_{p}^{j}}(C,H^{s}\centerdot R)\geq s$\\

or equivalently,\\

$I_{p_{j}}(C^{ns},\overline{H^{s}\centerdot R})=I_{p_{j}}(C^{ns},{\overline H}^{s}\centerdot\overline{R})\geq s$\\

for a suitable presentation $C^{ns}$ of a non-singular model of
$C$, see Lemma 5.12 of \cite{depiro1}, where we have used the
"lifted" form notation there. Using the method of conic
projections, see section 4 of \cite{depiro1}, we can find a plane
projective curve $C'$ birational to $C^{ns}$, such that the point
$p_{j}$ corresponds to a
non-singular point $q$ of $C'$ and;\\

$I_{p_{j}}(C^{ns},\overline{H}^{s}\centerdot\overline
R)=I_{q}(C',\overline{{\overline
H}^{s}\centerdot\overline{R}})=I_{q}(C',\overline{{\overline
H}}^{s}\centerdot\overline{\overline{R}})$\\

The result then follows by results of the paper \cite{depiro3} for
the intersections of plane projective curves. This shows $(**)$.
We now need to prove that, for an algebraic form $\phi_{\lambda}$
in $\Sigma$ and a branch $\gamma_{p}^{j}$ of $C$;\\

$I_{\gamma_{p}^{j}}(C,\phi_{\lambda}\centerdot
G)=I_{\gamma_{p}^{j}}(C,\phi_{\lambda})+I_{\gamma_{p}^{j}}(C,G)$\\

This follows by exactly the same argument, reducing to the case of
intersections between plane projective curves and using the
results of \cite{depiro3}. The result is then shown.

\end{proof}

\begin{theorem}{Birational Invariance of a $g_{n}^{r}$}\\

Let $\Phi:C_{1}\leftrightsquigarrow C_{2}$ be a birational map
between projective algebraic curves. Then, given a $g_{n}^{r}$ on
$C_{2}$, there exists a canonically defined $g_{n}^{r}$ on
$C_{1}$, depending only on the class $[\Phi]$ of the birational
map. Conversely, given a $g_{n}^{r}$ on $C_{1}$, there exists a
canonically defined $g_{n}^{r}$ on $C_{2}$, depending only on the
class $[\Phi^{-1}]$ of the birational map. Moreover, these
correspondences are inverse.

\end{theorem}

\begin{proof}
By Lemma 5.7 of \cite{depiro1}, $[\Phi]$ induces a bijection;\\

$[\Phi]^{*}:\bigcup_{O\in C_{2}}\gamma_{O}\rightarrow\bigcup_{O\in
C_{1}}\gamma_{O}$\\

of branches, with inverse given by ${[\Phi^{-1}]}^{*}$.

Then $[\Phi]^{*}$ extends naturally to a map on weighted sets of
degree $n$ by the formula;\\

$[\Phi]^{*}(n_{1}\gamma_{p_{1}}^{j_{1}}+\ldots+n_{r}\gamma_{p_{r}}^{j_{r}})=n_{1}[\Phi]^{*}(\gamma_{p_{1}}^{j_{1}})+\ldots+n_{r}[\Phi]^{*}(\gamma_{p_{r}}^{j_{r}})$\\

for a linear combination of branches
$\{\gamma_{p_{1}}^{j_{1}},\ldots,\gamma_{p_{r}}^{j_{r}}\}$ with\\
$n=n_{1}+\ldots+n_{r}$. Therefore, given a $g_{n}^{r}$ on $C_{2}$,
we obtain a canonically defined collection $[\Phi]^{*}(g_{n}^{r})$
of weighted sets on $C_{1}$ of degree $n$ $(*)$. It is trivial to
see that $[\Phi^{-1}]^{*}\circ[\Phi]^{*}(g_{n}^{r})$ recovers the
original $g_{n}^{r}$ on $C_{2}$, by the fact the map $[\Phi]^{*}$
on branches is invertible, with inverse given by
$[\Phi^{-1}]^{*}$. Let $C^{ns}$ be a non-singular model of $C_{1}$
and $C_{2}$ with morphisms $\Phi_{1}:C^{ns}\rightarrow C_{1}$ and
$\Phi_{2}:C^{ns}\rightarrow C_{2}$ such that
$\Phi\circ\Phi_{1}=\Phi_{2}$ and $\Phi^{-1}\circ\Phi_{2}=\Phi_{1}$
as birational maps (see the proof of Lemma 5.7 in \cite{depiro1}).
We then have that
$[\Phi]^{*}(g_{n}^{r})=[\Phi_{1}^{-1}]^{*}\circ[\Phi_{2}]^{*}(g_{n}^{r})$.
It remains to prove that this collection given by $(*)$ defines a
$g_{n}^{r}$ on $C_{1}$. We will prove first that
$[\Phi_{2}]^{*}(g_{n}^{r})$ defines a $g_{n}^{r}$ on $C^{ns}$
$(\dag)$. Let the original $g_{n}^{r}$ on $C_{2}$ be defined by a
linear system $\Sigma$, having finite intersection with $C_{2}$,
such that $dimension(\Sigma)=r$ and $degree(\Sigma)=n'$, after
removing some fixed branch contribution of multiplicity $n'-n$. We
may assume that $n'=n$, as if the fixed branch contribution in
question is given by $W_{0}$ and $g_{n}^{r}\cup W_{0}=g_{n'}^{r}$,
then $[\Phi_{2}]^{*}(g_{n}^{r})\cup
[\Phi_{2}]^{*}(W_{0})=[\Phi_{2}]^{*}(g_{n'}^{r})$, hence it is
sufficient to prove that $[\Phi_{2}]^{*}(g_{n'}^{r})$ defines a
$g_{n'}^{r}$. Let $W_{1}$ be the fixed branch contribution of the
$g_{n}^{r}$ on $C_{2}$ and let $g_{n''}^{r}\subset g_{n}^{r}$ be
obtained by removing this fixed branch contribution. It will be
sufficient to prove that $[\Phi_{2}]^{*}(g_{n''}^{r})$ defines a
$g_{n''}^{r}$ on $C^{ns}$ as
$[\Phi_{2}]^{*}(g_{n}^{r})=[\Phi_{2}]^{*}(g_{n''}^{r})\cup
[\Phi_{2}]^{*}(W_{1})$ and we may then use Lemma 3.13. Let
$\Phi_{\Sigma_{1}}$ and $\Phi_{\Sigma_{2}}$ be presentations of
the morphisms $\Phi_{1}$ and $\Phi_{2}$. We may assume that $
Base(\Sigma_{1})$ and $Base(\Sigma_{2})$ are disjoint. Let
$\{\overline{\phi_{\lambda}}\}$ denote the lifted family of forms
on $C^{ns}$, defined by the linear system $\Sigma$ and the
presentation $\Phi_{\Sigma_{2}}$. We claim that
$[\Phi_{2}]^{*}(g_{n''}^{r})$ is defined by this system after
removing its fixed branch contribution. In order to see this, we first show
that for any branch $\gamma_{p}^{j}$ of $C$;\\

$I_{\gamma_{p}^{j}}^{\Sigma,mobile}(C,\phi_{\lambda})=I_{p_{j}}^{\Sigma,mobile}(C^{ns},\overline{\phi_{\lambda}})$ $(*)$ $(1)$\\

where $p_{j}$ corresponds to $\gamma_{p}^{j}$ in the fibre
$\Gamma_{[\Phi_{2}]}(x,p)$, see Section 5 of \cite{depiro1}.
By Definition 2.20 and Lemma 5.23 of \cite{depiro1}, we have that;\\

$I_{p_{j}}^{\Sigma,mobile}(C^{ns},\overline{\phi_{\lambda}})=Card(C^{ns}\cap
({\mathcal V}_{p_{j}}\setminus
p_{j})\cap\overline{\phi_{\lambda'}})$ for
$\lambda'\in {\mathcal V}_{\lambda}$, generic\\
\indent \ \ \ \ \ \ \ \ \ \ \ \ \ \ \ \ \ \ \ \ \ \ \ \ \ \ \ \ \ \ \ \ \ \ \ \ \ \ \ \ \ \ \ \ \ \ \ \ \ \ \ \ \ \ \ \ \ \ \ \ \ \ \ \ in $Par_{\Sigma}$\\

$I_{\gamma_{p}^{j}}^{\Sigma,mobile}(C,\phi_{\lambda})=Card(C\cap
(\gamma_{p}^{j}\setminus p)\cap\phi_{\lambda'})$ for $\lambda'\in
{\mathcal V}_{\lambda}$, generic in $Par_{\Sigma}$\\

As $(\gamma_{p}^{j}\setminus p)$ is in biunivocal correspondence
with ${({\mathcal V}_{p_{j}}\setminus p_{j})}$ under the morphism
$\Phi_{2}$, we obtain immediately the result $(*)$. Now, using
Lemma 5.27 of \cite{depiro1}, we have that, if $\gamma_{p}^{j}$
appears in a weighted set $W_{\lambda}$ of the $g_{n''}^{r}$ with
multiplicity $s$, then the corresponding branch $\gamma_{p_{j}}$
appears in the weighted set $[\Phi_{2}]^{*}(W_{\lambda})$ with
multiplicity equal to
$s=I_{p_{j}}^{mobile}(C^{ns},\overline{\phi_{\lambda}})$. Again,
using Lemma 5.27 of \cite{depiro1}, we obtain that
$[\Phi_{2}]^{*}(W_{\lambda})$ is given by
$C^{ns}\sqcap\overline{\phi_{\lambda}}$, after removing all fixed
point contributions of the linear system $\Sigma$. We, therefore,
obtain that $[\Phi_{2}]^{*}(g_{n''}^{r})$ is defined by $\Sigma$,
after removing all fixed branch contributions, as required. This
proves $(\dag)$. We now claim that, for the given $g_{n}^{r}$ on
$C^{ns}$, $[\Phi_{1}^{-1}]^{*}(g_{n}^{r})$ defines a $g_{n}^{r}$ on
$C_{1}$, $(\dag\dag)$. Let $\Phi_{\Sigma_{3}}$ be a presentation
of the morphism $\Phi_{1}^{-1}$. If $\phi_{\lambda}$ is a form
belonging to the linear system $\Sigma$ defined on $C^{ns}$, using
the presentations $\Phi_{\Sigma_{1}}$ and $\Phi_{\Sigma_{3}}$ of
$\Phi_{1}$ and $\Phi_{1}^{-1}$, we obtain a lifted form
$\overline{\phi_{\lambda}}$ on $C_{1}$ and a lifted form
$\overline{\overline{\phi_{\lambda}}}$ on $C^{ns}$ again. We now
claim that, for $p\in C^{ns}$;\\

$I_{p}^{\Sigma,mobile}(C^{ns},\phi_{\lambda})=I_{p}^{\Sigma,mobile}(C^{ns},\overline{\overline{\phi_{\lambda}}})$ $(2)$\\

In order to see this, first observe that we can obtain the lifted
system of forms $\{\overline{\overline{\phi_{\lambda}}}\}$
directly from the linear system $\Sigma_{4}$, obtained by
composing bases of the linear systems $\Sigma_{1}$ and
$\Sigma_{3}$. The corresponding morphism $\Phi_{\Sigma_{4}}$
defines a birational map of $C^{ns}$ to itself, which is
equivalent to the identity map $Id$. Now the result follows
immediately from Definition 2.20 and Lemma 2.16 of \cite{depiro1},
both multiplicities are witnessed inside the canonical set $W$ of
$\Phi_{\Sigma_{4}}$, which, in this case, is just the domain of
definition of $\Phi_{\Sigma_{4}}$ on $C^{ns}$, see Definition 1.30
of \cite{depiro1}. Now, returning to the proof of $(\dag\dag)$, we
may suppose that the given $g_{n}^{r}$ on $C^{ns}$ is defined by
the linear system $\Sigma$, after removing all fixed branch
contributions. Combining $(1)$ and $(2)$, we have that, for a branch
 $\gamma_{p}^{j}$ of $C_{1}$;\\

$I_{\gamma_{p}^{j}}^{\Sigma,mobile}(C_{1},\overline{\phi_{\lambda}})=I_{p_{j}}^{\Sigma,mobile}(C^{ns},\overline{\overline{\phi_{\lambda}}})=I_{p_{j}}^{\Sigma,mobile}(C^{ns},\phi_{\lambda})$\\

The result now follows from the same argument as above, using
Lemma 5.27 of \cite{depiro1}. This completes the theorem.

\begin{rmk}
Using the quoted Theorem 1.33 of \cite{depiro1}, one can use the
Theorem to reduce calculations involving $g_{n}^{r}$ on projective
algebraic curves to calculations on plane projective curves. This idea
is central to the philosophy of the "Italian School" of algebraic geometry.
\end{rmk}
\end{proof}

We finally note the following;\\

\begin{lemma}
For a given $g_{n}^{r}$, we always have that $r\leq n$.
\end{lemma}

\begin{proof}
The proof is almost identical to Lemma 2.24 of \cite{depiro1}. We
leave the details to the reader.

\end{proof}

\end{section}
\begin{section}{A Theory of Complete Linear Series on an Algebraic
Curve}

We now develop further the theory of $g_{n}^{r}$ on an algebraic
curve $C$, analogously to classical results for divisors on
non-singular algebraic curves. We will first assume that $C$ is a
plane projective algebraic curve, defined by some homogeneous
polynomial $F(X,Y,Z)$. Without loss of generality, we will use the
coordinates $x=X/Z$ and $y=Y/Z$ for local calculations on the
curve $C$, defined in this system by $f(x,y)=0$. Using Theorem
3.14, we will later derive general results for $g_{n}^{r}$
on an algebraic curve from the corresponding calculations for the plane case.\\

We consider first the case when $r=1$. By results of the previous
section, a $g_{n}^{1}$ is defined by a pencil $\Sigma$ of
algebraic curves $\{\phi(x,y)+\lambda\phi'(x,y)=0\}_{\lambda\in
P^{1}}$ (in affine coordinates), after removing some fixed branch
contribution, where, by convention, we interpret the algebraic
curve $\phi(x,y)+\infty\phi'(x,y)=0$ to be $\phi'(x,y)=0$. We
assume that the $g_{n}^{1}$ is, in fact, cut out by this pencil.
Now suppose that $\gamma_{p}$ is a branch of $C$. We may assume
that $p$ corresponds to the origin $O$ of the affine coordinate
system $(x,y)$, (use a linear transformation and the result of
Lemma 4.1) By Theorem 6.1 of \cite{depiro1}, we can find algebraic
power series $\{x(t),y(t)\}$, with $x(t)=y(t)=0$, parametrising
$\gamma_{p}$. We can now substitute the power series in order to
obtain a formal expression of the form;\\

${\phi(x(t),y(t))\over\phi'(x(t),y(t))}={t^{i}u(t)\over
t^{j}v(t)}=t^{i-j}u(t)v(t)^{-1}$, where
$\{u(t),v(t),u(t)v(t)^{-1}\}$\\
\indent \ \ \ \ \ \ \ \ \ \ \ \ \ \ \ \ \ \ \ \ \ \ \ \ \ \ \ \ \
\ \ \ \ \ \ \ \ \ \ \ \ \ \ \ \ \ \ \ \  are units in $L[[t]]$.

We then define;\\

$(i)$. \ \ $ord_{\gamma_{p}}({\phi\over\phi'})=i-j$,\\

\indent \ \ \ \ \ \ \ \ $val_{\gamma_{p}}({\phi\over\phi'})=0$,\ \ \ \ \ \ \ \ \ \ \ \ \ \ \ \ \ \ \ if $i>j$, (${\phi\over\phi'}$ has a zero of order $i-j$)\\

$(ii)$. \ $ord_{\gamma_{p}}({\phi\over\phi'})=j-i$,\\

\indent \ \ \ \ \ \ \ \ \ \ $val_{\gamma_{p}}({\phi\over\phi'})=\infty$,\ \ \ \ \ \ \ \ \ \ \ \ \ \ \ \ \ \ if $i<j$, (${\phi\over\phi'}$ has a pole of order $j-i$)\\

$(iii)$. $ord_{\gamma_{p}}({\phi\over\phi'})=ord_{t}(h(t)-h(0))$,\\

\indent \ \ \ \ \ \ \ \ $val_{\gamma_{p}}({\phi\over\phi'})=h(0)$,\ \ \ \ \ \ \ \ \ \ \ \ \ if $i=j$ and $h(t)=u(t)v(t)^{-1}$\\

Observe that in all cases, $ord_{\gamma_{p}}$ gives a
\emph{positive} integer, while $val_{\gamma_{p}}$ determines an
element of $P^{1}$. In order to see that this construction does
not depend on the particular power series representation of
the branch, we require the following lemma;\\

\begin{lemma}
Let $\{C,\gamma_{p},\phi,\phi',g_{n}^{1},\Sigma\}$ be as defined above, then;\\

$ord_{\gamma_{p}}({\phi\over\phi'})=I_{\gamma_{p}}(C,\phi-\lambda\phi')$,
\ \ \ \ \ \ if $\gamma_{p}$ is not a base branch for the
$g_{n}^{1}$ \indent \ \ \ \ \ \ \ \ \ \ \ \ \ \ \ \ \ \ \ \ \ \ \
\ \ \ \ \ \
\ \ \ \ \ \ \ \ \ \ \ and ${\phi\over\phi'}(p)=val_{\gamma_{p}}({\phi\over\phi'})=\lambda$.\\

$ord_{\gamma_{p}}({\phi\over\phi'})=I_{\gamma_{p}}^{\Sigma,mobile}(C,\phi-\lambda\phi')$,
if $\gamma_{p}$ is a base branch for the $g_{n}^{1}$ and\\
\indent \ \ \ \ \ \ \ \ \ \ \ \ \ \ \ \ \ \ \  \ \ \ \ \ \ \  \ \
\  \ \ \ \ \ \ \ \ \ \ \ $\lambda=val_{\gamma_{p}}({\phi\over\phi'})$ is unique such that,\\
\indent\ \ \ \ \ \ \ \ \ \ \ \ \ \ \ \ \ \ \ \ \ \ \ \ \ \ \ \ \ \ \ \ \ \ \ \ \ \ \ \  for $\mu\neq\lambda$;\\
\indent \ \ \ \ \ \ \ \ \ \ \ \ \ \ \ \ \ \ \ \ \ \ \ \ \ \ \ \ \
 \ \ \ \ \ \ \ \ \ \ \
$I_{\gamma_{p}}(C,\phi-\lambda\phi')>I_{\gamma_{p}}(C,\phi-\mu\phi').$

\end{lemma}

\begin{proof}
Suppose that $\gamma_{p}$ is not a base branch for the
$g_{n}^{1}$, then ${\phi\over\phi'}(p)=\lambda$ is well defined,
if we interpret $(c/0)=\infty$ for $c\neq 0$, and
$\phi-\lambda\phi'$ is the unique curve in the pencil passing
through $p$. It is trivial to check, using the facts that
$\phi(p)=\phi(x(0),y(0))$ and $\phi'(p)=\phi'(x(0),y(0))$, that,
in all cases, $val_{\gamma_{p}}({\phi\over\phi'})=\lambda$ as
well. By Theorem 6.1 of \cite{depiro1}, we have that;\\

$I_{\gamma_{p}}(C,\phi-\lambda\phi')=ord_{t}[(\phi-\lambda\phi')(x(t),y(t))]$\\

If $\lambda=0$, then $\phi(p)=0$ and $\phi'(p)\neq 0$, hence, by a
straightforward algebraic calculation,
$\phi(x(t),y(t))=t^{i}u(t)$, for some $i\geq 1$, and
$\phi'(x(t),y(t))=v(t)$ for $\{u(t),v(t)\}$ units in $L[[t]]$.
Therefore,
$ord_{\gamma_{p}}({\phi\over\phi'})=ord_{t}\phi(x(t),y(t))$ and
the result follows.\\

If $\lambda=\infty$, then $\phi(p)\neq 0$ and $\phi(p)=0$, hence,
$\phi(x(t),y(t))=u(t)$ and $\phi'(x(t),y(t))=t^{j}v(t)$, for some
$j\geq 1$, and $\{u(t),v(t)\}$ units in $L[[t]]$. Therefore,
$ord_{\gamma_{p}}({\phi\over\phi'})=ord_{t}\phi'(x(t),y(t))$ and
the result follows.\\

If $\lambda\neq\{0,\infty\}$, then $\phi(x(t),y(t))=u(t)$ and
$\phi'(x(t),y(t))=v(t)$ with $\{u(t),v(t)\}$ units in $L[[t]]$. As
$v(t)$ is a unit in $L[[t]]$, we have that;\\

$ord_{t}({u(t)\over v(t)}-{u(0)\over
v(0)})=ord_{t}(v(t)({u(t)\over
v(t)}-{u(0)\over v(0)}))=ord_{t}(u(t)-{u(0)\over v(0)}v(t))$\\

Hence, by definition of $ord_{\gamma_{p}}$;\\

$ord_{\gamma_{p}}({\phi\over\phi'})=ord_{t}[(\phi-\lambda\phi')(x(t),y(t))]$\\

and the result follows.\\

Now suppose that $\gamma_{p}$ is a base branch for the
$g_{n}^{1}$, then $\phi(p)=\phi'(p)=0$ and we have that
$\phi(x(t),y(t))=t^{i}u(t)$ and $\phi'(x(t),y(t))=t^{j}v(t)$, for
some $i,j\geq 1$ and $\{u(t),v(t)\}$ units in $L[[t]]$. Again, we
divide the proof into the following cases;\\

$i>j$. In this case, by definition, $val_{\gamma_{p}}({\phi\over\phi'})=0$. We compute;\\

$ord_{t}(\phi(x(t),y(t))-\lambda\phi'(x(t),y(t)))=ord_{t}(t^{i}u(t)-\lambda
t^{j}v(t))$\\

When $\lambda=0$, we obtain, by Theorem 6.1 of \cite{depiro1},
that $I_{\gamma_{p}}(C,\phi)=i$ and, for $\lambda\neq 0$, that
$I_{\gamma_{p}}(C,\phi-\lambda\phi')=j$. Using Lemma 5.27 of
\cite{depiro1}, we obtain that
$I_{\gamma_{p}}^{\Sigma,mobile}(C,\phi)=i-j=ord_{\gamma_{p}}({\phi\over\phi'})$, as required.\\

$i<j$. In this case, by definition,
$val_{\gamma_{p}}({\phi\over\phi'})=\infty$. The computation for
$ord_{\gamma_{p}}$ is similar, with the critical value
being $\lambda=\infty$.\\

$i=j$. We compute;\\

$ord_{t}(\phi(x(t),y(t))-\lambda\phi'(x(t),y(t)))=ord_{t}[t^{i}(u(t)-\lambda
v(t))]$\\

Again, there exists a unique value of $\lambda={u(0)\over
v(0)}=val_{\gamma_{p}}({\phi\over\phi'})\neq\{0,\infty\}$ such
that $ord_{t}(u(t)-\lambda v(t))=k\geq 1$. By the same calculation
as above, we have that
$I_{\gamma_{p}}^{\Sigma,mobile}(C,\phi-\lambda\phi')=k$, for this
critical value of $\lambda$. By a similar algebraic calculation to
the above, using the fact that $v(t)$ is a unit, we also compute
$ord_{\gamma_{p}}({\phi\over\phi'})=k$, hence the result
follows.\\
\end{proof}

We now show the following;\\

\begin{lemma}
Given any algebraic curve $C\subset P^{w}$, with function field
$L(C)$, for a non-constant rational function $f\in L(C)$ and a
branch $\gamma_{p}$, we can unambiguously define
$ord_{\gamma_{p}}(f)$ and $val_{\gamma_{p}}(f)$.
\end{lemma}

\begin{proof}
The proof is similar to the above. We may, without loss of
generality, assume that $p$ corresponds to the origin of a
coordinate system $(x_{1},\ldots,x_{w})$. Using Theorem 6.1 of
\cite{depiro1}, we can find algebraic power series
$(x_{1}(t),\ldots,x_{w}(t))$ parametrising the branch
$\gamma_{p}$. By the assumption that $f$ is non-constant, we can
find a representation of $f$ as a rational function
$\phi(x_{1},\ldots,x_{w})\over\phi'(x_{1},\ldots,x_{w})$ in this
coordinate system, such that the pencil $\Sigma$ defined by
$\{\phi,\phi'\}$ has finite intersection with $C$, hence defines a
$g_{n}^{1}$. Using the method above, we can define
$ord_{\gamma_{p}}({\phi\over\phi'})$ and
$val_{\gamma_{p}}({\phi\over\phi'})$ for this representation. The
proof of Lemma 4.1 shows that these are defined independently of
the particular power series parametrising the branch. We need to
check that they are also defined independently of the particular
representation of $f$. Suppose that
$\{\phi_{1},\phi_{2},\phi_{3},\phi_{4}\}$ are algebraic forms with
the property that
${\phi_{1}\over\phi_{2}}={\phi_{3}\over\phi_{4}}$ as rational
functions on $C$. We claim that, for any branch $\gamma_{p}$ of
$C$,
$ord_{\gamma_{p}}({\phi_{1}\over\phi_{2}})=ord_{\gamma_{p}}({\phi_{3}\over\phi_{4}})$
and
$val_{\gamma_{p}}({\phi_{1}\over\phi_{2}})=val_{\gamma_{p}}({\phi_{3}\over\phi_{4}})$,
$(*)$. In order to see this, let $U\subset NonSing(C)$ be an open
subset of $C$, on which ${\phi_{1}\over\phi_{2}}$ and
${\phi_{3}\over\phi_{4}}$ are defined and equal. Let $g_{n}^{1}$
and $g_{m}^{1}$ on $C$ be defined by the pencils
$\Sigma_{1}=\{\phi_{1}-\lambda\phi_{2}\}_{\lambda\in P^{1}}$ and
$\Sigma_{2}=\{\phi_{3}-\lambda\phi_{4}\}_{\lambda\in P^{1}}$. Let
$V=U\setminus Base(\Sigma_{1})\cup Base(\Sigma_{2})$. Then
$V\subset U$ is also an open subset of $C$, which we will refer to
as the canonical set. Now, suppose that $\gamma_{p}\subset V$. We
will prove $(*)$ for this branch. As both
${\phi_{1}\over\phi_{2}}$ and ${\phi_{3}\over\phi_{4}}$ are
defined and equal at $p$, using the argument in Lemma 4.1, we have
that
$val_{\gamma_{p}}({\phi_{1}\over\phi_{2}})=val_{\gamma_{p}}({\phi_{3}\over\phi_{4}})$.
It is therefore sufficient, again by Lemma
4.1, to show that;\\

$I_{\gamma_{p}}(C,\phi_{1}-\lambda\phi_{2})=I_{\gamma_{p}}(C,\phi_{3}-\lambda\phi_{4})$,
for ${\phi_{1}\over\phi_{2}}(p)={\phi_{3}\over\phi_{4}}(p)=\lambda$ $(\dag)$\\

Suppose that $I_{\gamma_{p}}(\phi_{1}-\lambda\phi_{2})=m$, then,
by Lemma 5.25 of \cite{depiro1}, we can find $\lambda'\in{\mathcal
V}_{\lambda}\cap P^{1}$ and $\{p_{1},\ldots,p_{m}\}=V\cap
{\mathcal V}_{p}\cap (\phi_{1}-\lambda'\phi_{2})=0$ witnessing
this multiplicity. As $\{p,p_{1},\ldots,p_{m}\}$ lie inside $V$,
we also have that $\{p_{1},\ldots,p_{m}\}\subset V\cap{\mathcal
V}_{p}\cap(\phi_{3}-\lambda'\phi_{4})=0$, hence
$I_{\gamma_{p}}(C,\phi_{3}-\lambda\phi_{4})\geq m$. The result
$(\dag)$ then follows from the converse argument.\\

Now, suppose that $\gamma_{p}$ is one of the finitely many
branches of $C$, not lying inside $V$. We will just consider the
case when $\gamma_{p}$ is a base branch for both the $g_{n}^{1}$
and the $g_{m}^{1}$ defined above, the other cases being similar.
In order to prove $(*)$ for this branch, it is sufficient, by
Lemma 4.1, to show that;\\

$I_{\gamma_{p}}^{\Sigma_{1},mobile}(C,\phi_{1}-\lambda\phi_{2})=I_{\gamma_{p}}^{\Sigma_{2},mobile}(C,\phi_{3}-\mu\phi_{4})$,
for the critical values\\
\indent \ \ \ \ \ \ \ \ \ \ \ \ \ \ \ \ \ \ \ \ \ \ \ \ \ \ \ \ \ \ \ \ \ \ \ \ \ \ \ \ \ \ \ \ \ \ \ \ \ \ \ \ \ \ \ \ \ \ \ \ \ \ $\{\lambda,\mu\}$\\

and that the critical values $\{\lambda,\mu\}$ coincide,
$(\dag\dag)$.\\

Using the argument to prove $(\dag)$, witnessing the corresponding
multiplicities in the canonical set $V$, it follows that for
$\emph{any}$ $\nu\in P^{1}$;\\

$I_{\gamma_{p}}^{\Sigma_{1},mobile}(C,\phi_{1}-\nu\phi_{2})=I_{\gamma_{p}}^{\Sigma_{2},mobile}(C,\phi_{3}-\nu\phi_{4})$, $(\dag\dag\dag)$\\

If the critical values $\{\lambda,\mu\}$ were distinct, we would
have that;\\

$I_{\gamma_{p}}^{\Sigma_{1},mobile}(C,\phi_{1}-\lambda\phi_{2})>I_{\gamma_{p}}^{\Sigma_{1},mobile}(C,\phi_{1}-\mu\phi_{2})$\\
\indent\ \ \ \ \ \ \ \ \ \ \ \ \ \ \ \ \ $||$\ \ \ \ \ \ \ \ \ \ \
\ \ \ \ \ \ \ \ \ \ \ \ \ \ \ \ \ \ \ \ $||$\ \ \ \ \ \ \ \ \
\\
\indent $I_{\gamma_{p}}^{\Sigma_{2},mobile}(C,\phi_{3}-\lambda\phi_{4})<I_{\gamma_{p}}^{\Sigma_{2},mobile}(C,\phi_{3}-\mu\phi_{4})$\\

which is clearly a contradiction. Hence, $\lambda=\mu$ and the
result $(\dag\dag)$ follows from $(\dag\dag\dag)$. The lemma is
shown.
\end{proof}

\begin{lemma}{Birational Invariance of $ord_{\gamma_{p}}$ and $val_{\gamma_{p}}$}\\

Let $\Phi:C_{1}\leftrightsquigarrow C_{2}$ be a birational map
between projective algebraic curves with corresponding
isomorphisms $\Phi^{*}:L(C_{2})\rightarrow L(C_{1})$ and\\
$[\Phi]^{*}:\bigcup_{p\in C_{2}}\gamma_{p}\rightarrow\bigcup_{q\in
C_{1}}\gamma_{q}$ . Then, for non-constant $f\in L(C_{2})$ and
$\gamma_{p}$ a branch of $C_{2}$,
$ord_{\gamma_{p}}(f)=ord_{[\Phi]^{*}\gamma_{p}}(\Phi^{*}f)$ and
$val_{\gamma_{p}}(f)=val_{[\Phi]^{*}\gamma_{p}}(\Phi^{*}f)$.

\end{lemma}

\begin{proof}
Let $f$ be represented as a rational function by
${\phi_{1}\over\phi_{2}}$, as in Lemma 4.2, and consider the
$g^{1}_{n}$ on $C_{2}$, defined by the linear system
$\Sigma=\{\phi_{1}-\lambda\phi_{2}\}_{\lambda\in P^{1}}$. Let
$\Phi_{\Sigma_{1}}$ be a presentation of the birational map
$\Phi$. Using this presentation, we may lift the system $\Sigma$
to a corresponding linear system
$\{\overline{\phi_{1}}-\lambda\overline{\phi_{2}}\}_{\lambda\in
P^{1}}$. It is trivial to check that $\Phi^{*}f$ is represented by
the rational function
${\overline{\phi_{1}}\over\overline{\phi_{2}}}$. The proof of
Theorem 3.14 shows that, for a branch $\gamma_{p}$ of $C_{2}$;\\

$I_{\gamma_{p}}^{\Sigma,mobile}(C_{2},\phi_{1}-\lambda\phi_{2})=I_{[\Phi]^{*}\gamma_{p}}^{\Sigma,mobile}(C_{1},\overline{\phi_{1}}-\lambda\overline{\phi_{2}})$, $(*)$\\

We now need to consider the following cases;\\

Case 1. $\gamma_{p}$ and $[\Phi]^{*}\gamma_{p}$ are not base
branches for $\Sigma$ on $C_{2}$ and $C_{1}$.\\

Case 2. $\gamma_{p}$ is not a base branch, but
$[\Phi]^{*}\gamma_{p}$ is a base branch for $\Sigma$ on\\
\indent \ \ \ \ \ \ \ \ \ \ \  $C_{2}$ and $C_{1}$.\\

Case 3. $\gamma_{p}$ is a base branch and $[\Phi]^{*}\gamma_{p}$
is a base branch for $\Sigma$ on $C_{2}$\\
\indent\ \ \ \ \ \ \ \ \ \ \  and $C_{1}$.\\

For Case 1, we have, by Lemma 4.1 and $(*)$;\\

$ord_{\gamma_{p}}({\phi_{1}\over\phi_{2}})=I_{\gamma_{p}}(C_{2},\phi_{1}-\lambda\phi_{2})=I_{[\Phi]^{*}\gamma_{p}}(C_{1},\overline{\phi_{1}}-\lambda\overline{\phi_{2}})=ord_{[\Phi]^{*}\gamma_{p}}({\overline{\phi_{1}}\over\overline{\phi_{2}}})$\\

where
${\phi_{1}\over\phi_{2}}(p)={\overline{\phi_{1}}\over\overline{\phi_{2}}}(q)=val_{\gamma_{p}}({\phi_{1}\over\phi_{2}})=val_{\gamma_{q}}({\overline{\phi_{1}}\over\overline{\phi_{2}}})=\lambda$
and $[\Phi]^{*}\gamma_{p}=\gamma_{q}$.\\

For Case 3, we have, by Lemma 4.1, $(*)$ and a similar argument to
the previous lemma, in order to show the critical value
$\lambda=val_{\gamma_{p}}({\phi_{1}\over\phi_{2}})$ is also the
critical value
$val_{\gamma_{q}}({\overline{\phi_{1}}\over\overline{\phi_{2}}})$
for the lifted system at the
corresponding branch $[\Phi]^{*}\gamma_{p}$, that;\\

$ord_{\gamma_{p}}({\phi_{1}\over\phi_{2}})=I_{\gamma_{p}}^{\Sigma,mobile}(C_{2},\phi_{1}-\lambda\phi_{2})=I_{[\Phi]^{*}\gamma_{p}}^{\Sigma,mobile}(C_{1},\overline{\phi_{1}}-\lambda\overline{\phi_{2}})=ord_{[\Phi]^{*}\gamma_{p}}({\overline{\phi_{1}}\over\overline{\phi_{2}}})$\\

Case 2 is similar, we leave the details to the reader.\\

The lemma now follows from the previous lemma, that the
definitions of $ord_{\gamma_{p}}(f)$,
$ord_{[\Phi]^{*}\gamma_{p}}(\Phi^{*}f)$,$val_{\gamma_{p}}(f)$ and
$val_{[\Phi]^{*}\gamma_{p}}(\Phi^{*}f)$ are independent of their
particular representations.

\end{proof}

We now show;\\

\begin{lemma}{flatness}\\

Let $C$ be a projective algebraic curve, then, to any non-constant
rational function $f$ on $C$, we can associate a $g_{n}^{1}$ on
$C$, which we will denote by $(f)$, where $n=deg(f)$.
\end{lemma}

\begin{proof}
We define the weighted set $(f=\lambda)$ as follows;\\

$(f=\lambda):=\{n_{\gamma_{1}},\ldots,n_{\gamma_{r}}\}$\\

where
$\{\gamma_{1},\ldots,\gamma_{r}\}=\{\gamma:val_{\gamma}(f)=\lambda\}$ and $n_{\gamma}=ord_{\gamma}(f)$.\\

As $\lambda$ varies over $P^{1}$, we obtain a series of weighted
sets $W_{\lambda}$ on $C$. We claim that this series does in fact
define a $g_{n}^{1}$. In order to see this, let $f$ be represented
as a rational function by $\phi\over\phi'$. As before, we consider
the pencil $\Sigma$ of forms defined by
$(\phi-\lambda\phi')_{\lambda\in P^{1}}$. We claim that the series
is defined by this system $\Sigma$, after removing its fixed
branch contribution, $(*)$. In order to see this, we compare the
weighted sets $(f=\lambda)$ and $C\sqcap (\phi-\lambda\phi')$. For
a branch $\gamma_{p}$ which is not a fixed branch of the system
$\Sigma$,
we have, using Lemmas 4.1 and 4.2, that;\\

 $\gamma_{p}\in (f=\lambda)$ iff $val_{\gamma_{p}}(f)=\lambda$ iff
${\phi\over\phi'}(p)=\lambda$ iff $p\in
C\cap(\phi-\lambda\phi')$\\

In this case, by Lemmas 4.1 and 4.2, we have that;\\

$n_{\gamma_{p}}=ord_{\gamma_{p}}(f)=ord_{\gamma_{p}}({\phi\over\phi'})=I_{\gamma_{p}}(C,\phi-\lambda\phi')$\\

For a branch $\gamma_{p}$ which is a fixed branch of the system
$\Sigma$, we have, by Lemmas 4.1 and 4.2, that;\\

$\gamma_{p}\in (f=\lambda)$ iff
$val_{\gamma_{p}}({\phi\over\phi'})=\lambda$ iff $p\in
C\cap(\phi-\lambda\phi')$ and $\lambda$ is a critical value for
the system $\Sigma$ at $\gamma_{p}$.\\

In this case, by Lemmas 4.1 and 4.2, we have that;\\

$n_{\gamma_{p}}=ord_{\gamma_{p}}(f)=ord_{\gamma_{p}}({\phi\over\phi'})=I_{\gamma_{p}}^{\Sigma,mobile}(C,\phi-\lambda\phi')$ $(1)$\\

Let $I_{\gamma_{p}}=min_{\mu\in
P^{1}}I_{\gamma_{p}}(C,\phi-\mu\phi')$ be the fixed branch
contribution of $\Sigma$ at $\gamma_{p}$. Then, at the critical
value $\lambda$ for the system $\Sigma$;\\

 $I_{\gamma_{p}}^{\Sigma,mobile}(\phi-\lambda\phi')=I_{\gamma_{p}}(C,\phi-\lambda\phi')-I_{\gamma_{p}}$ $(2)$\\

Hence, the result $(*)$ follows from $(1),(2)$ and the definition
of\\
 $C\sqcap(\phi-\lambda\phi')$.\\

Finally, we show that $n=deg(f)$. Let $\Gamma_{f}$ be the
correspondence determined by the rational map $f:C\rightsquigarrow
P^{1}$. By classical arguments, $deg(f)$ is equal to the
cardinality of the generic fibre $\Gamma_{f}(\lambda)$, for
$\lambda\in P^{1}$. Fixing a presentation ${\phi\over\phi'}$ for
$f$, if $U\subset NonSing(C)$ is the canonical set for this
presentation, one may assume that the generic fibre
$\Gamma_{f}(\lambda)$ lies inside $U$. By Lemma 2.17 of
\cite{depiro1}, one may also assume that the corresponding
weighted set of the $g_{n}^{1}$ defined by $(f=\lambda)$ consists
of $n$ distinct branches, centred at the points of the generic
fibre $\Gamma_{f}(\lambda)$. Therefore, the result follows.

\end{proof}

\begin{rmk}
By convention, for a non-zero rational function $c\in
{L\setminus\{0\}}$, we define $(c=0)$ and $(c=\infty)$ to be the
empty weighted sets. The notion of a weighted set in a
$g_{n}^{1}$, generalises the classical notion of the divisor on a
non-singular curve. Using the above theorem, we can make sense of
the notion of linear equivalence of weighted sets.
\end{rmk}

We make the following definition;\\

\begin{defn}{Linear equivalence of weighted sets}\\

Let $C$ be an algebraic curve and let $A$ and $B$ be weighted sets
on $C$ of the same total multiplicity. We define $A\equiv B$ if
there exists a $g_{n}^{r}$ on $C$ such that $A$ and $B$ belong to
this $g_{n}^{r}$ as weighted sets.

\end{defn}

\begin{theorem}
Let hypotheses be as in the previous definition. If $A\equiv B$,
then there exists a rational function $g$ on $C$, such that $A$ is
defined by $(g=0)$ and $B$ is defined by $(g=\infty)$, possibly
after adding some fixed branch contribution.

\end{theorem}

\begin{proof}
If $r=0$ in the definition, then we must have that $A=B$. Hence,
we obtain the statement of the theorem by adding the fixed branch
contribution $A$ to the empty $g_{0}^{0}$, defined by
$(c=0)=(c=\infty)$, for a non-constant $c\in L^{*}$. Otherwise, by
the definition of a $g_{n}^{r}$, we may, without loss of
generality, find a pencil $\Sigma$ of algebraic forms,
$\{\phi-\lambda\phi'\}_{\lambda\in P^{1}}$, having finite intersection with $C$, such that;\\

$A=C\sqcap(\phi-\lambda_{1}\phi)$,\\

$B=C\sqcap(\phi-\lambda_{2}\phi')$\ \ \ \ \  $(\lambda_{1}\neq
\lambda_{2})$\\

Let $f$ be the rational function on $C$ defined by
${\phi\over\phi'}$. If $A$ and $B$ have no branches in common
(with multiplicity), $(\dag)$, then the pencil $\Sigma$ can have
no fixed
branches and, by Lemma 4.4, we have that;\\

$A=(f=\lambda_{1})$\\

$B=(f=\lambda_{2})$\ \ \ \ \ \ \ $(\lambda_{1}\neq\lambda_{2})$\\

Now we can find an algebraic automorphism $\alpha$ of $P^{1}$,
taking $\lambda_{1}$ to $0$ and $\lambda_{2}$ to $\infty$. We will
assume that $\{\lambda_{1},\lambda_{2}\}\neq\infty$, in which case
$\alpha$ can be given, for a coordinate $z$ on $P^{1}$, by the
Mobius transformation ${z-\lambda_{1}\over z-\lambda_{2}}$. The
other cases are left to the reader. Let $g$ be the rational
function on $C$ defined by $\alpha\circ f$. Now, suppose that
$\gamma$ is a branch of $C$, with $val_{\gamma}(f)=\lambda$ and
$ord_{\gamma}(f)=m$. Then, we claim that
$val_{\gamma}(g)=\alpha(\lambda)$ and $ord_{\gamma}(g)=m$, $(*)$.
If $\lambda\neq\{\lambda_{2},\infty\}$, using the method before
Lemma 4.1, we obtain
the following power series representation of $g$ at $\gamma$;\\

${(\lambda+\mu t^{m}+o(t^{m}))-\lambda_{1}\over(\lambda+\mu
t^{m}+o(t^{m}))-\lambda_{2}}=[(\lambda-\lambda_{1})+\mu t^{m}+o(t^{m})]\centerdot{1\over (\lambda-\lambda_{2})}[1-{\mu\over(\lambda-\lambda_{2})}t^{m}+o(t^{m})]$\\
\indent \ \ \ \ \ \ \ \ \ \ \ \ \ \ \ \ \ \ \ \ \
$={\lambda-\lambda_{1}\over
\lambda-\lambda_{2}}+t^{m}[{\mu(\lambda-\lambda_{2})-\mu(\lambda-\lambda_{1})\over
(\lambda-\lambda_{2})^{2}}]+o(t^{m})$\\
\indent \ \ \ \ \ \ \ \ \ \ \ \ \ \ \ \ \ \ \ \ \
$={\lambda-\lambda_{1}\over
\lambda-\lambda_{2}}+t^{m}[{\mu(\lambda_{1}-\lambda_{2})\over
(\lambda-\lambda_{2})^{2}}]+o(t^{m})$\\

and the claim $(*)$ follows from the assumption that
$\lambda_{1}\neq \lambda_{2}$. If $\lambda=\lambda_{2}$, we obtain
the following power series representation of $g$ at
$\gamma$;\\

${(\lambda+\mu t^{m}+o(t^{m}))-\lambda_{1}\over(\mu
t^{m}+o(t^{m}))}={1\over t^{m}}\centerdot
[(\lambda-\lambda_{1})+\mu t^{m}+o(t^{m})]\centerdot [\mu+o(1)]^{-1}$\\

which gives that $val_{\gamma}(g)=\infty=\alpha(\lambda_{2})$ and
$ord_{\gamma}(g)=m$, using the fact that $\lambda\neq\lambda_{1}$.
Finally, if $\lambda=\infty$, the Mobius transformation at
$\infty$ is given by ${{1\over z}-\lambda_{1}\over {1\over
z}-\lambda_{2}}={1-\lambda_{1}z\over 1-\lambda_{2}z}$ and $g$ may
be represented at $\gamma$ by
${\phi-\lambda_{1}\phi'\over\phi-\lambda_{2}\phi'}$. We then
obtain the power series representation of $g$ at $\gamma$;\\

${(t^{i}u(t)-\lambda_{1}t^{i+m}v(t))\over(t^{i}u(t)-\lambda_{2}t^{i+m}v(t))}={(u(t)-\lambda_{1}t^{m}v(t))\over(u(t)-\lambda_{2}t^{m}v(t))}={[1-\lambda_{1}t^{m}{v(t)\over
u(t)}]\over [1-\lambda_{2}t^{m}{v(t)\over u(t)}]}$\\
$\indent \ \ \ \ \ \ \ \ \ \ \ \ \ \ \ \ \ \ \ \ =1+(\lambda_{2}-\lambda_{1})t^{m}w(t)+o(t^{m})$, for $\{u(t),v(t),w(t)\}$\\
\indent \ \ \ \ \ \ \ \ \ \ \ \ \ \ \ \ \ \ \ \ \ \ \ \ \ \ \ \ \ \ \ \ \ \ \ \ \ \ \ \ \ \ \ \ \ \ \ \ \ \ \ \ \ \ \ \ \ \ \ \ \ \ \ \ \ \ \ \  units in $L[[t]]$\\

which gives that $val_{\gamma}(g)=1=\alpha(\infty)$ and
$ord_{\gamma}(g)=m$, using the fact that
$\lambda_{1}\neq\lambda_{2}$ again. This gives the claim $(*)$. It
follows that the weighted sets $(f=\lambda)$ correspond exactly to
the weighted sets $(g=\alpha(\lambda))$, in particularly the
$g_{n}^{1}$ defined by $(f)$ and $(g)$, as in Lemma 4.4, is the
same. With this new parametrisation of the $g_{n}^{1}$, we then
have that;\\

$A=(g=0)$\\

$B=(g=\infty)$\\

Hence, the result follows, with the assumption $(\dag)$. If $A$
and $B$ have branches in common, with multiplicity, we let $A\cap
B$ denote the weighted set consisting of these common branches
(with multiplicity). Then, the same argument holds, replacing $A$
by $A\setminus B=A-(A\cap B)$ and $B$ by $B\setminus A=B-(A\cap
B)$. After adding the fixed branch contribution $(A\cap B)$ to the
$g_{n}^{1}$ defined by $(g)$, we then obtain the result. Note
that, by Lemma 3.13, this addition defines a $g_{n+n'}^{1}$, where
$n'$ is the total multiplicity of $(A\cap B)$. \\

\end{proof}

\begin{rmk}
The definition we have given of linear equivalence of weighted
sets on a projective algebraic curve $C$ generalises the modern
definition of linear equivalence for effective divisors on a
smooth projective algebraic curve. More precisely we have;\\

Modern Definition; Let $A$ and $B$ be effective divisors on a
smooth projective algebraic curve $C$, then $A\equiv B$ iff
$A-B=div(g)$, for some $g\in L(C)^{*}$.\\

See, for example, p161 of \cite{Shaf} for relevant definitions and
notation. We now show that our definition is the same in this
case. First, observe that there exists a natural bijection between
the set of effective divisors on $C$, in the sense of \cite{Shaf},
and the collection of weighted sets on $C$, $(*)$. This follows
immediately from the fact, given in Lemma 5.29 of \cite{depiro1},
that, for each point $p\in C$, there exists a unique branch
$\gamma_{p}$, centred at $p$. Secondly, observe that the notion of
$div(g)$, for $g\in L(C)$, as given in \cite{Shaf}, is the same as
the notion of $div(g)$ which we give in Definition 4.9 below,
(taking into account the identification $(*)$), $(\dag)$. This
amounts to checking that, for a point $p\in C$, with corresponding
branch
$\gamma_{p}$;\\

$v_{p}(g)=ord_{\gamma_{p}}(g)$ $(\dag\dag)$\\

where $v_{p}(g)$ is defined in p152 of \cite{Shaf}. First, one can
use the fact, given in Lemma 4.9 of \cite{depiro1}, together with
remarks from the final section of this paper, that there exists a
birational map $\phi:C\leftrightsquigarrow C'$, such that $C'$ is
a plane projective algebraic curve, and $p$ corresponds to a
non-singular point $p'\in C'$ with $\{p,p'\}$ lying inside the
canonical sets associated to $\phi$. Using the calculation given
below, in Lemma 4.10, for $ord_{\gamma_{p}}$, and the definition
of $v_{p}$, one can assume that $v_{p}(g)\geq 0$ and $g\in
O_{p,C}$. Let $g'\in L(C')$ denote the corresponding rational
function to $g$ on $L(C)$. It is then a trivial algebraic
calculation, using the fact that the local rings $O_{p,C}$ and
$O_{p',C'}$ are isomorphic, to show that $v_{p}(g)=v_{p'}(g')$. It
also follows from Lemma 4.3 that
$ord_{\gamma_{p}}(g)=ord_{\gamma_{p'}}(g')$. Hence, it is
sufficient to check $(\dag\dag)$ for the plane projective curve $C'$.
We may, without loss of generality, assume that $v_{p'}(g')\geq 1$
and that $g'$ is represented in some choice of affine coordinates
$\{x,y\}$ by the polynomial $q(x,y)$. If $Q(X,Y,Z)$ denotes the
projective equation of this polynomial and $p'$ corresponds to the
origin of this coordinate system, then;\\

$v_{p'}(g')=I_{p'}(C,Q)=length({L[x,y]\over <h,q>})$\\

where $h$ is a defining equation for $C'$ in the coordinate system
$\{x,y\}$ and $I_{p'}$ is the algebraic intersection multiplicity.
It also follows from Lemma 4.1, that;\\

$ord_{\gamma_{p'}}(g')=I_{\gamma_{p'}}(C,Q)$\\

Hence, it is sufficient to check that;\\

$I_{p'}(C,Q)=I_{\gamma_{p'}}(C,Q)$\\

This calculation was done in Theorem 2.10, hence
$(\dag\dag)$ and therefore $(\dag)$ is shown. Thirdly, it remains
to check that the definitions of linear equivalence are the same.
In order to see this, observe that we can write (for effective
divisors or weighted sets $A$ and $B$);\\

$A-B=(A\setminus B)+(A\cap B)]-[(B\setminus A)+(A\cap B)]=(A\setminus B)-(B\setminus A)$, $(\dag\dag\dag)$\\

If $A\equiv B$ in the sense of weighted sets (Definition 4.6),
then the calculation $(\dag\dag\dag)$ (which removes the fixed
branch contribution) and Theorem 4.7 shows that $A-B=div(g)$, for
some rational function $g\in L(C)$, where, here, $div(g)$ is as
defined in Definition 4.9. By $(\dag)$, it then follows that
$A\equiv B$ as effective divisors. Conversely, if $A\equiv B$ as
effective divisors, then there exists a rational function $g\in
L(C)$ such that $A-B=div(g)$, in the sense of the modern
definition given above. The above calculations $(\dag\dag\dag)$
and $(\dag)$ then show that $div(g)=(A\setminus B)-(B\setminus
A)$, in the sense of Definition 4.9 below. It follows, by Lemma
4.4, that there exists a $g_{n}^{1}$ to which $(A\setminus B)$ and
$(B\setminus A)$ belong as weighted sets. Adding the fixed branch
contribution $(A\cap B)$ to this $g_{n}^{1}$, we then obtain that
$A\equiv B$ in the sense of Definition 4.6, as required.
\end{rmk}

\begin{defn}
Let $C$ be a projective algebraic curve and let $f$ be a non-zero
rational function on $C$. Then we define $div(f)$ to
be the weighted set $A-B$ where;\\

$A=(f=0)$,\indent $B=(f=\infty)$\\

\end{defn}

We now require the following lemma;\\

\begin{lemma}
Let $C$ be a projective algebraic curve, and let $f$ and $g$ be
non-zero rational functions on $C$. Then;\\

$div({1\over f})=-div(f)$\\

$div(fg)=div(f)+div(g)$\\

$div({f\over g})=div(f)-div(g)$\\

\end{lemma}

\begin{proof}
In order to prove the first claim, it is sufficient to show that,
for a branch $\gamma$ of $C$;\\

$val_{\gamma}(f)=0$ iff $val_{\gamma}({1\over f})=\infty$\\

$val_{\gamma}(f)=\infty$ iff $val_{\gamma}({1\over f})=0$\\

and $ord_{\gamma}$ is preserved in both cases. This follows
trivially from the relevant power series calculation at a branch.
Namely, we can represent $f$ by ${\phi\over\phi'}$ and ${1\over
f}$ by ${\phi'\over\phi}$. Substituting the branch
parametrisation, we obtain that;\\

$val_{\gamma}(f)=0, ord_{\gamma}(f)=m$ iff $f\sim t^{m}u(t)$, \
$m\geq 1,u(t)\in L[[t]]$ a
unit.\\
\indent \ \ \ \ \ \ \ \ \ \ \ \ \ \ \ \ \ \ \ \ \ \ \ \ \ \ \ \ \
\ \ \ \ \ iff ${1\over f}\sim
t^{-m}u(t)^{-1}$\\
\indent \ \ \ \ \ \ \ \ \ \ \ \ \ \ \ \ \ \ \ \ \ \ \ \ \ \ \ \ \
\ \ \ \ \ iff $val_{\gamma}(f)=\infty, ord_{\gamma}(f)=m$\\

and the calculation for $val_{\gamma}(f)=\infty,
ord_{\gamma}(f)=m$ is similar.\\

In order to prove the second claim, we need to verify the following cases at a branch $\gamma$ of $C$;\\

Case 1. If $val_{\gamma}(f)=val_{\gamma}(g)\in\{0,\infty\}$,
$ord_{\gamma}(f)=m$ and
$ord_{\gamma}(g)=n$\\

\indent \ \ \ \ \ \ \ \ \ \ \  then $val_{\gamma}(fg)\in\{0,\infty\}$ and $ord_{\gamma}(fg)=m+n$\\

Case 2. If $val_{\gamma}(f)\neq val_{\gamma}(g)\in\{0,\infty\}$,
$ord_{\gamma}(f)=m$ and $ord_{\gamma}(g)=n$\\

\indent \ \ \ \ \ \ \ \ \ \ \ \ then
$val_{\gamma}(fg)\in\{0,\infty\}$ and $ord_{\gamma}(fg)=|m-n|$\\

Case 3. If exactly one of $val_{\gamma}(f)$ and $val_{\gamma}(g)$
is in $\{0,\infty\}$, with\\
\indent \ \ \ \ \ \ \  \ \ \ \ $ord_{\gamma}(f)$ or $ord_{\gamma}(g)=m$\\

\indent \ \ \ \ \ \ \ \ \ \ \ then $val_{\gamma}(fg)\in\{0,\infty\}$, with $ord_{\gamma}(fg)=m$.\\

Case 4. If neither of $val_{\gamma}(f)$ and $val_{\gamma}(g)$ are
in $\{0,\infty\}$\\

\indent \ \ \ \ \ \ \ \ \ \ \ \ then $val_{\gamma}(fg)$ is not in
$\{0,\infty\}$\\

If $f$ is represented by ${\phi\over\phi'}$ and $g$ is represented
by ${\psi\over\psi'}$, then we can represent $fg$ by
${\phi\psi\over\phi'\psi'}$. The proof of these cases then follow by elementary power series calculations at the branch $\gamma$.
For example, for Case 2, if $val_{\gamma}(f)=0$ and $ord_{\gamma}(f)=m$, $val_{\gamma}(g)=\infty$ and $ord_{\gamma}(g)=n$, then we have;\\

$f\sim t^{n}u(t)$, $g\sim t^{-m}v(t)$, $fg\sim
t^{n}t^{-m}u(t)v(t)=t^{n-m}w(t)$,\\
\indent \ \ \ \ \ \ \ \ \ \ \ \ \ \ \ \ \ \ \ \ \ \ \ \ \ \ \ \ \
\ \ \ \ \ \ \ \  for $\{u(t),v(t),w(t)\}$ units in
$L[[t]]$.\\

The third claim follows from the first two claims.

\end{proof}

We now claim the following;\\

\begin{theorem}{Transitivity of Linear Equivalence}\\

Let $C'$ be an algebraic curve. If $A,B,C$ are weighted sets on
$C'$ of the same total multiplicity, then, if $A\equiv B$ and
$B\equiv C$, we must have that $A\equiv C$.

\end{theorem}

\begin{proof}
By Theorem 4.7, we can find rational functions $f$ and $g$ on
$C'$, such that;\\

$(A\setminus B)-(B\setminus A)=div(f)$\\

$(B\setminus C)-(C\setminus B)=div(g)$\\

By Lemma 4.10, we have that;\\

$div(fg)=(A\setminus B)-(B\setminus A)+(B\setminus C)-(C\setminus
B)$\\

By drawing a Venn diagram, one easily checks that;\\

$(A\setminus B)-(B\setminus A)=(A\cap B^{c}\cap C^{c})+(A\cap
B^{c}\cap C)-(A^{c}\cap B\cap C^{c})-\\
\indent \ \ \ \ \ \ \ \ \ \ \ \ \ \ \ \ \ \ \ \ \ \ \ \ \ \ \ (A^{c}\cap B\cap C)$\\
\indent \ \ \ \ \ \ \ \ \ \ $+$\\
\indent $(B\setminus C)-(C\setminus B)=(A\cap B\cap
C^{c})+(A^{c}\cap
B\cap C^{c})-(A^{c}\cap B^{c}\cap C)-\\
\indent \ \ \ \ \ \ \ \ \ \ \ \ \ \ \ \ \ \ \ \ \ \ \ \ \ \ \
(A\cap B^{c}\cap C)$\\
\indent \ \ \ \ \ \ \ \ \ \ \ $||$\\
\indent $(A\setminus C)-(C\setminus A)=(A\cap B^{c}\cap
C^{c})+(A\cap
B\cap C^{c})-(A^{c}\cap B^{c}\cap C)-\\
\indent \ \ \ \ \ \ \ \ \ \ \ \ \ \ \ \ \ \ \ \ \ \ \ \ \ \ \ (A^{c}\cap B\cap C)$\\

Hence, $div(fg)=(A\setminus C)-(C\setminus A)$. Now, given the
$g_{n}^{1}$ defined by the rational function $fg$, as in Lemma
4.4, it follows that $(A\setminus C)$ and $(C\setminus A)$ belong
to this $g_{n}^{1}$ as weighted sets. We can now add the fixed
branch contribution $A\cap C$ to this $g_{n}^{1}$, giving a
$g_{n+n'}^{1}$, to which $A$ and $C$ belong as weighted sets.
Therefore, the result follows.

\end{proof}

As an immediate corollary, we have;\\

\begin{theorem}
Let $C$ be a projective algebraic curve, then $\equiv$ is an
equivalence relation on weighted sets for $C$ of a given
multiplicity.

\end{theorem}

We also have;\\

\begin{theorem}{Linear Equivalence preserved by Addition}\\

Let $C'$ be a projective algebraic curve and suppose that
$\{A,B,C,D\}$ are weighted sets on $C'$ with;\\

$A\equiv B$ and $C\equiv D$\\

then;\\

$A+C\equiv B+D$\\

\end{theorem}

\begin{proof}
By Definition 4.6, we can find a $g_{n}^{r}$ containing $C$ and
$D$ as weighted sets. If $s$ is the total multiplicity of $A$,
then, by Lemma 3.13, we can add the weighted set $A$ as a fixed
branch contribution to this $g_{n}^{r}$ and obtain a
$g_{n+s}^{r}$, containing $A+C$ and $A+D$ as weighted sets. Hence,
by Definition 4.6 again, we have that;\\

$A+C\equiv A+D$ $(1)$\\

Similarily, one shows, by adding $D$ as a fixed branch contribution to the $g_{n'}^{r'}$ containing
$A$ and $B$ as weighted sets, that;\\

$A+D\equiv B+D$ $(2)$\\

The result then follows immediately by combining $(1)$, $(2)$ and
using Theorem 4.11.

\end{proof}

We now develop further the theory of $g_{n}^{r}$ on a projective
algebraic curve $C$. We begin with the following definition;\\

\begin{defn}{Subordinate $g_{n}^{r}$}\\

Let  $\{g_{n}^{r},g_{n}^{t}\}$ be given on $C$ with the
\emph{same} order $n$. Then we say that;\\

$g_{n}^{r}\subseteq g_{n}^{t}$\\

if \emph{every} weighted set in $g_{n}^{r}$ is included in the
weighted sets of the $g_{n}^{t}$.

\end{defn}

We now claim the following;\\

\begin{theorem}{Amalgamation of $g_{n}^{r}$}\\

Let $\{g_{n}^{r},g_{n}^{s}\}$ be given on $C$, having a common
weighted set $G$, then there exists $t$ with $r\leq t, s\leq t$
and a $g_{n}^{t}$ such that $g_{n}^{r}\subseteq g_{n}^{t}$ and
$g_{n}^{s}\subseteq g_{n}^{t}$.
\end{theorem}

\begin{proof}
Assume first that $\{g_{n}^{r},g_{n}^{s}\}$ have no fixed branch
contribution and are defined exactly by linear systems. Then we
can find algebraic forms
$\{\phi_{0},\psi_{0}\}$ such that;\\

$G=(C\sqcap\phi_{0}=0)=(C\sqcap\psi_{0}=0)$\\

and;\\

$g_{n}^{r}$ is defined by
$C\sqcap(\epsilon_{0}\phi_{0}+\epsilon_{1}\phi_{1}+\ldots+\epsilon_{r}\phi_{r}=0)$\\

$g_{n}^{s}$ is defined by $C\sqcap(\eta_{0}\psi_{0}+\eta_{1}\psi_{1}+\ldots+\eta_{s}\psi_{s}=0)$\\

Now consider the linear system $\Sigma$ defined by;\\

$\epsilon\phi_{0}\psi_{0}+\psi_{0}(\epsilon_{1}\phi_{1}+\ldots+\epsilon_{r}\phi_{r})+\phi_{0}(\eta_{1}\psi_{1}+\ldots+\eta_{s}\psi_{s})=0$\\

and let $g_{m}^{t}$ be defined by $\Sigma$. As
$deg(\psi_{0}\phi_{0})=deg(\psi_{0})+deg(\phi_{0})$, we have that
$m=2n$. We claim that the fixed branch contribution of
$g_{2n}^{t}$ is exactly $G$, $(*)$. In order to see this, observe
that we can write an algebraic form in $\Sigma$ as;\\

$\psi_{0}\phi_{\bar\epsilon}+\phi_{0}\psi_{\bar\eta}$\\

If $\gamma$ is a branch counted $w$-times in $G$,
then, using the proof at the end of Lemma 3.13 and linearity of multiplicity at a branch,
see \cite{depiro1};\\

$I_{\gamma}(C,\psi_{0}\phi_{\bar\epsilon})=I_{\gamma}(C,\psi_{0})+I_{\gamma}(C,\phi_{\bar\epsilon})\geq
w$\\

$I_{\gamma}(C,\phi_{0}\psi_{\bar\eta})=I_{\gamma}(C,\phi_{0})+I_{\gamma}(C,\psi_{\bar\eta})\geq
w$\\

$I_{\gamma}(C,\psi_{0}\phi_{\bar\epsilon}+\phi_{0}\psi_{\bar\eta})=
min\{I_{\gamma}(C,\psi_{0}\phi_{\bar\epsilon}),I_{\gamma}(C,\phi_{0}\psi_{\bar\eta})\}\geq
w$ $(\dag)$\\

Hence, $\gamma$ is $w$-fold for the $g_{2n}^{t}$ and $G$ is
contained in the fixed branch contribution of the $g_{2n}^{t}$. In
order to obtain the exactness statement, $(*)$, first observe
that, if $\gamma$ is a fixed branch of the $g_{2n}^{t}$, then, in
particular, it belongs to $(C\sqcap \phi_{0}\psi_{0}=0)$. Hence,
it belongs either to $(C\sqcap \phi_{0}=0)$ or
$(C\sqcap\psi_{0}=0)$. Hence, it belongs to $G$. Now, using the
fact that the original $\{g_{n}^{r},g_{n}^{s}\}$ had no fixed
branch contribution, we can easily find $\phi_{\bar\epsilon_{0}}$
and $\psi_{\bar\eta_{0}}$ with $G$ disjoint from both $(C\sqcap
\phi_{\bar\epsilon_{0}}=0)$ and $(C\sqcap\psi_{\bar\eta_{0}}=0)$.
Then, by the same argument $(\dag)$, we obtain, for a branch $\gamma$ of $G$;\\

$I_{\gamma}(C,\psi_{0}\phi_{\bar\epsilon_{0}}+\phi_{0}\psi_{\bar\eta_{0}})=w$\\

hence, $\gamma$ is counted $w$-times in
$C\sqcap(\psi_{0}\phi_{\bar\epsilon_{0}}+\phi_{0}\psi_{\bar\eta_{0}}=0)$
and, therefore, $(*)$ holds, as required. Now, as $G$ had total
multiplicity $n$, removing this fixed branch contribution from the
$g_{2n}^{t}$, we obtain a $g_{n}^{t}$. We then claim that
$g_{n}^{r}\subseteq g_{n}^{t}$ and $g_{n}^{s}\subseteq g_{n}^{t}$,
$(**)$. By Definition 4.14, it is sufficient to check that, if
$\{W_{1},W_{2}\}$ are weighted sets appearing in
$\{g_{n}^{r},g_{n}^{s}\}$, defined by $(C\sqcap
\phi_{\bar\epsilon}=0)$ and $(C\sqcap\psi_{\bar\eta}=0)$, then
they appear in the $g_{n}^{t}$. We clearly have that both
$\psi_{0}\phi_{\bar\epsilon}$ and $\phi_{0}\psi_{\bar\eta}$ belong
to $\Sigma$ and the calculation $(\dag)$ shows that;\\

$C\sqcap(\psi_{0}\phi_{\bar\epsilon}=0)=W_{1}+G$\\

$C\sqcap(\phi_{0}\psi_{\bar\eta}=0)=W_{2}+G$\\

Hence, the result $(**)$ follows after removing the fixing branch
contribution $G$. The fact that $r\leq t$ and $s\leq t$ then
follows easily from the definition of the dimension of a $g_{n}^{r}$
and Theorem 3.3.\\

Now consider the case when the $\{g_{n}^{r},g_{n}^{s}\}$ are
defined exactly by linear systems and \emph{have} a fixed branch
contribution. Let $G_{1}\subseteq G$ and $G_{2}\subseteq G$ be
these fixed branch contributions and let $G_{3}=G_{1}\cap G_{2}$.
We claim that the fixed branch contribution of the $g_{2n}^{t}$
defined by $\Sigma$, as given above, in this case is exactly
$G_{3}+G$. The proof is similar to the above and left to the
reader. Now, removing the fixed branch contribution $G$, we obtain
a series $g_{n}^{t}$ with fixed branch contribution $G_{3}$. A
similar proof to the above, left to the reader, shows that this
$g_{n}^{t}$ contains the original series
$\{g_{n}^{r},g_{n}^{s}\}$. Finally, we need to consider the case
when the $\{g_{n}^{r},g_{n}^{s}\}$ are defined, after removing
some fixed branch contribution from linear series. Let $G_{1}$ and
$G_{2}$, with total multiplicity $r_{1}$ and $r_{2}$, be these
fixed branch contributions and let
$\{g_{n+r_{1}}^{r},g_{n+r_{2}}^{s}\}$ be the series obtained from
adding these fixed branch contributions to
$\{g_{n}^{r},g_{n}^{s}\}$. In this case, the linear system
$\Sigma$, as given above, defines a $g_{2n+r_{1}+r_{2}}^{t}$. We
claim that the weighted set $G\cup G_{1}\cup G_{2}$, of total
multiplicity $(n+r_{1}+r_{2})$, is contained in the fixed branch
contribution of this series. This follows from a similar
calculation, using the method above, the details are left to the
reader. Removing this weighted set from the
$g_{2n+r_{1}+r_{2}}^{t}$, we obtain a $g_{n}^{t}$ and a similar
calculation shows that this contains the original
$\{g_{n}^{r},g_{n}^{s}\}$, again the details are left to the
reader.
\end{proof}

As a corollary, we have;\\

\begin{theorem}
Let a $g_{n}^{r}$ be given on $C$, then there exists a
\emph{unique} $g_{n}^{t}$ on $C$, with $r\leq t\leq n$, such that;\\

$g_{n}^{r}\subseteq g_{n}^{t}$\\

and, for \emph{any} $g_{n}^{s}$ such that $g_{n}^{r}\subseteq
g_{n}^{s}$, we have that;\\

$g_{n}^{s}\subseteq g_{n}^{t}$\\
\end{theorem}

\begin{proof}
By Lemma 3.16, we can find $r\leq t\leq n$ and a $g_{n}^{t}$ on
$C$, with $g_{n}^{r}\subseteq g_{n}^{t}$ and $t$ maximal with this
property. If $g_{n}^{r}\subseteq g_{n}^{s}$, then
$\{g_{n}^{s},g_{n}^{t}\}$ would contain a common weighted set. By
Theorem 4.15, we could then find $t'\leq n$ such that $s\leq t'$,
$t\leq t'$ and $g_{n}^{s}\subseteq g_{n}^{t'}$,
$g_{n}^{t}\subseteq g_{n}^{t'}$. If $g_{n}^{s}\varsubsetneq
g_{n}^{t}$, then, by elementary dimension considerations, we would
have that $t<t'\leq n$ and $g_{n}^{r}\subset g_{n}^{t'}$,
contradicting maximality of $t$. Hence, $g_{n}^{s}\subseteq
g_{n}^{t}$. The uniqueness statement also follows from a similar
amalgamation argument, using Theorem 4.15.

\end{proof}

We can then make the following definition;\\

\begin{defn}
We call a $g_{n}^{r}$ on $C$ complete if it cannot be strictly
contained in a $g_{n}^{t}$ of greater dimension. If $G$ is any
weighted set on $C$ of total multiplicity $n$, then we define
$|G|$ to be the unique complete $g_{n}^{t}$ to which $G$ belongs.
\end{defn}

We then have that;\\

\begin{theorem}
Let $G$ be a weighted set on $C$, then, $G\equiv G'$ if and only
if $G'$ belongs to $|G|$. In particular, $G\equiv G'$ if and only
if $|G|=|G'|$.

\end{theorem}

\begin{proof}
The proof of the first part of the theorem is quite
straightforward. By definition, if $G'$ belongs to $|G|$, then
$G\equiv G'$. Conversely, if $G'\equiv G$, then, by Definition
4.6, we can find a $g_{n}^{1}$, containing the given weighted sets
$G$ and $G'$. By Theorem 4.16, we can find a unique complete
$g_{n}^{t}$ on $C$, with $1\leq t\leq n$, such that
$g_{n}^{1}\subseteq g_{n}^{t}$. As $G$ belongs to this $g_{n}^{t}$
as a weighted set, it follows by Definition 4.17 that
$|G|=g_{n}^{t}$. Hence, $G'$ belongs to $|G|$ as required. For the
second part, if $G\equiv G'$, then, by the first part, $G'$
belongs to $|G|$. It follows immediately from Definition 4.17 and
Theorem 4.16, that $|G|\subseteq |G'|$. Reversing this argument,
we have that $|G'|\subseteq |G|$, hence $|G|=|G'|$ as required.
Conversely, if $|G|=|G'|$, then clearly $G\equiv G'$ by Definition
4.6.
\end{proof}

We now make the following definition;\\

\begin{defn}{Linear System of a Weighted Set}\\

Let $G$ be a weighted set on a projective algebraic curve $C$,
then we define the Riemann-Roch space ${\mathcal L}(C,G)$ or
${\mathcal L}(G)$ to be the vector space defined as;\\

$\{g\in L(C)^{*}:div(g)+G\geq 0\}\cup\{0\}$\\

where $div(g)$ was defined in Definition 4.9.

\end{defn}

\begin{rmk}
That ${\mathcal L}(G)$ defines a vector space follows easily from
Lemma 4.10, the fact that, for non-constant rational functions
$\{f,g,f+g\}\subset
L(C)$ and a branch $\gamma$ of $C$, we have that;\\

$ord_{\gamma}(f+g)\geq min\{ord_{\gamma}(f),ord_{\gamma}(g)\}$, $(*)$\\

where, for this remark only, $ord_{\gamma}$ is counted
\emph{negatively} if $val_{\gamma}$ is infinite, and an argument
on constants, $(**)$.
 We now give a brief proof of $(*)$;\\

We just consider the following $2$ cases;\\

Case 1. $val_{\gamma}(f)<\infty$ and $val_{\gamma}(g)<\infty$\\

We then have, substituting the relative parametrisations, that;\\

$f\sim c+c_{1}t^{m}+\ldots$ and $g\sim d+d_{1}t^{n}+\ldots$, where
$ord_{\gamma}(f)=m\geq 1$, $ord_{\gamma}(g)=n\geq 1$ and
$\{c_{1},d_{1}\}\subset L$ are non-zero. Then;\\

$f+g\sim (c+d)+c_{1}t^{m}+d_{1}t^{n}+\ldots$\\

If $(f+g)-(c+d)\equiv 0$, as an algebraic power series in
$L[[t]]$, then $(f+g)=(c+d)$ as a rational function on $C$,
contradicting the assumption. Hence, we obtain that
$ord_{\gamma}(f+g)=min\{ord_{\gamma}(f),ord_{\gamma}(g)\}$, if
$m\neq n$ or $m=n$ and $c_{1}+d_{1}\neq 0$, and
$ord_{\gamma}(f+g)>min\{ord_{\gamma}(f),ord_{\gamma}(g)\}$
otherwise. Hence, $(*)$ is shown in this case.\\

Case 2. $val_{\gamma}(f)=val_{\gamma}(g)=\infty$\\

We then have that;\\

$f\sim c_{1}t^{-m}+\ldots$ and $g\sim d_{1}t^{-n}+\ldots$, where
$ord_{\gamma}(f)=-m\leq -1$, $ord_{\gamma}(g)=-n\leq -1$ and
$\{c_{1},d_{1}\}\subset L$ are non-zero. Then;\\

$f+g\sim c_{1}t^{-m}+d_{1}t^{-n}+\ldots$\\

By the assumption that $f+g$ is not a constant, if $m=n$ and
$c_{1}+d_{1}=0$, we must have higher order terms in $t$ in the
Cauchy series for $(f+g)$, hence
$ord_{\gamma}(f+g)>min\{ord_{\gamma}(f),ord_{\gamma}(g)\}$.
Otherwise, we have that
$ord_{\gamma}(f+g)=min\{ord_{\gamma}(f),ord_{\gamma}(g)\}$, hence
$(*)$ is shown in this case as well.\\

The remaining cases are left to the reader. One should also
consider the case of constants, $(**)$. Technically, one cannot
define $ord_{\gamma}$ for a constant in $L$. However, we did, by
convention, define $div(c)=0$, for $c\in L^{*}$, in Remarks 4.5.\\
\end{rmk}

We now show the following;\\

\begin{lemma}
For a weighted set $G$, $dim({\mathcal L}(G))=t+1$, where $t$ is
given in Definition 4.17. In particular, ${\mathcal L}(G)$ is
finite dimensional.

\end{lemma}

\begin{proof}
Let $t$ be given by Definition 4.17. If $t=0$, then $G=(0)$ and
${\mathcal L}(G)=L$. This follows easily from the well known fact
that the only regular functions on a projective algebraic curve
are the constants (see, for example, \cite{Shaf}, p59). In this
case, we then have that $dim({\mathcal L}(G))=1$, as required.
Otherwise, let $t\geq 1$ be given as in Definition 4.17, with the
unique complete $g_{n}^{t}$ containing $G$. After adding some
fixed branch contribution $W$, we can find a linear system
$\Sigma$, having finite intersection with $C$, with basis
$\{\phi_{0},\ldots,\phi_{j},\ldots,\phi_{t}\}$ defining this
$g_{n}^{t}$. Moreover, we may assume that $C\sqcap \phi_{0}=G\cup
W$, $(*)$. Let $\{f_{1},\ldots,f_{j},\ldots,f_{t}\}$ be the
sequence of rational functions on $C$ defined by
$f_{j}={\phi_{j}\over\phi_{0}}$. We claim that;\\

$div(f_{j})+G\geq 0$, for $1\leq j\leq t$ $(**)$\\

In order to show $(**)$, it is sufficient to prove that, for a
branch $\gamma$ with $val_{\gamma}(f_{j})=\infty$, we have that
$\gamma$ belong to $G$ and, moreover, that $\gamma$ is counted at
least $ord_{\gamma}(f_{j})$ times in $G$. Let $\Sigma_{j}$ be the
pencil of forms defined by $(\phi_{j}-\lambda\phi_{0})_{\lambda\in
P^{1}}$. By the proof of Lemma 4.4, we have that $(f_{j}=\infty)$
is defined by $(C\sqcap\phi_{0})$, after removing the fixed branch
contribution of this pencil. By $(*)$ and the fact that the fixed
branch contribution of $\Sigma_{j}$ includes $W$, we have that
$(f_{j}=\infty)\subseteq G$. Hence, $(**)$ is shown as required.
By Definition 4.19, we then have that $f_{j}$ belongs to
${\mathcal L}(G)$. We now claim that there do \emph{not} exist
constants $\{c_{0},\ldots,c_{j},\ldots,c_{t}\}\subset L$ such that;\\

$c_{0}+c_{1}f_{1}+\ldots+c_{j}f_{j}+\ldots+c_{t}f_{t}=0$ $(***)$\\

as rational functions on $C$. If so, we would have that;\\

$c_{0}\phi_{0}+c_{1}\phi_{1}+\ldots+c_{j}\phi_{j}+\ldots+c_{t}\phi_{t}$\\

vanished identically on $C$, contradicting the fact that $\Sigma$
has finite intersection with $C$. Hence, by $(***)$,
$\{1,f_{1},\ldots,f_{t}\}\subset {\mathcal L}(G)$ are linearly
independent and $dim({\mathcal L}(G))\geq t+1$. Conversely,
suppose that $dim({\mathcal L}(G))\geq k+1$, then we can find
$\{1,f_{1},\ldots,f_{j},\ldots,f_{k}\}\subset {\mathcal L}(G)$
which are linearly independent, $(\dag)$. By the usual method of
equating denominators, we can find algebraic forms
$\{\phi_{0},\ldots,\phi_{k}\}$ of the same degree, such that
$f_{j}$ is represented by $\phi_{j}\over\phi_{0}$, for $1\leq
j\leq k$. Let $\Sigma$ be the linear system defined by this
sequence of forms. By $(\dag)$, $\Sigma$ has finite intersection
with $C$. Let $W$, having total multiplicity $n'$, be the fixed
branch contribution of this system and let
$(C\sqcap\phi_{0})=G_{0}\cup W$. We claim that $G_{0}\subseteq G$,
$(\dag\dag)$. Suppose not, then there exists a branch $\gamma$
with $I_{\gamma}^{\Sigma,mobile}(C,\phi_{0})=s$, where $\gamma$ is
counted strictly less than $s$-times in $G$. By the definition of
$I_{\gamma}^{\Sigma,mobile}$, we can find a form $\phi_{\lambda}$
belonging to $\Sigma$, distinct from $\phi_{0}$, witnessing this
multiplicity. Consider the pencil $\Sigma_{\lambda}$ defined by
$(\phi_{\lambda}-\mu\phi_{0})_{\mu\in P^{1}}$. We then clearly
have that $I_{\gamma}^{\Sigma_{\lambda},mobile}(C,\phi_{0})=s$ as
well, $(\dag\dag\dag)$. Let
$f_{\lambda}={\phi_{\lambda}\over\phi_{0}}$. By the proof of Lemma
4.4, we have that $(f_{\lambda}=\infty)$ is defined by
$(C\sqcap\phi_{0})$, after removing the fixed branch contribution
of $\Sigma_{\lambda}$. By $(\dag\dag\dag)$, it follows that the
branch $\gamma$ is counted $s$-times in $(f_{\lambda}=\infty)$ and
therefore $div(f_{\lambda})+G\ngeq 0$. However, $f_{\lambda}$ is a
linear combination of $\{1,\ldots,f_{k}\}$, hence
$f_{\lambda}\in{\mathcal L}(G)$, which is a contradiction. Hence,
$(\dag\dag)$ is shown. Now, consider the $g_{n}^{k}$ defined by
$\Sigma$. Let $W'$ be the weighted set $G\setminus G_{0}$ of total
multiplicity $n''$. By Lemma 3.13, we can add the weighted set
$W'$ to the $g_{n}^{k}$ and obtain a $g_{n+n''}^{k}$ with fixed
branch contribution $W'\cup W$. Now, removing the fixed branch
contribution $W$ from this $g_{n+n'}^{k}$, we obtain a
$g_{n+n''-n'}^{k}$ containing $G$ exactly as a weighted set. It
follows, from Definition 4.17, that $k\leq t$. Hence, in
particular, $dim({\mathcal L}(G))$ is finite and $dim({\mathcal
L}(G)\leq t+1$. Therefore, the lemma is proved.

\end{proof}

We now extend the notion of linear equivalence to include virtual,
or non-effecive, weighted sets.\\

\begin{defn}
We define a generalised weighted set $G$ on $C$ to be a linear
combination of branches;\\

$n_{1}\gamma_{p_{1}}^{j_{1}}+\ldots+n_{r}\gamma_{p_{r}}^{j_{r}}$\\

where $\{n_{1},\ldots,n_{r}\}$ belong to ${\mathcal Z}$. If
$\{n_{1},\ldots,n_{r}\}$ belong to ${\mathcal Z}_{\geq 0}$, we
call the weighted set effective. Otherwise, we call the weighted
set virtual. We define $n=n_{1}+\ldots+n_{r}$ to be the total
multiplicity or degree of $G$.
\end{defn}

\begin{rmk}
It is an easy exercise to see that there exist well defined
operations of addition and subtraction on generalised weighted
sets. It is also easy to check that any generalised weighted set
$G$ may be written uniquely as $G_{1}-G_{2}$, where
$\{G_{1},G_{2}\}$ are \emph{disjoint effective} weighted sets.
\end{rmk}

\begin{defn}
Let $A$ and $B$ be generalised weighted sets on $C$ of the same
total multiplicity. Let $\{A_{1},A_{2}\}$ and $\{B_{1},B_{2}\}$ be
the unique effective weighted sets, as given by the previous
remark. Then we define;\\

$(A_{1}-A_{2})\equiv (B_{1}-B_{2})$ iff $(A_{1}+B_{2})\equiv
(B_{1}+A_{2})$\\

and;\\

$A\equiv B$ iff $(A_{1}-A_{2})\equiv (B_{1}-B_{2})$\\

\end{defn}

\begin{rmk}
Note that if $\{A_{1}',A_{2}'\}$ and $\{B_{1}',B_{2}'\}$ are
\emph{any} effective weighted sets such that;\\

$A=A_{1}'-A_{2}'$ and $B=B_{1}'-B_{2}'$\\

then $A\equiv B$ iff $A_{1}'+B_{2}'\equiv B_{1}'+A_{2}'$\\

The proof is just manipulation of effective weighted sets. We
clearly have that;\\

$A_{1}+A_{2}'=A_{1}'+A_{2}$ and $B_{1}+B_{2}'=B_{1}'+B_{2}$ $(*)$\\

We then have;\\

\indent $A\equiv B$\indent iff \ \ \ \ \ \ \ $A_{1}+B_{2}\equiv B_{1}+A_{2}$(Definition 4.24)\\
\indent \ \ \ \ \ \ \ \ \ \ \ iff $A_{1}+A_{2}'+B_{2}\equiv
B_{1}+A_{2}+A_{2}'$
(Theorem 4.13)\\
\indent \ \ \ \ \ \ \ \ \ \ \ iff $A_{1}'+A_{2}+B_{2}\equiv
B_{1}+A_{2}+A_{2}'$ (by (*))\\
\indent \ \ \ \ \ \ \ \ \ \ \ iff $\ \ \ \ \ \ \
A_{1}'+B_{2}\equiv B_{1}+A_{2}'$
(Theorem 4.13)\\
\indent \ \ \ \ \ \ \ \ \ \ \ iff $A_{1}'+B_{2}+B_{1}'\equiv
B_{1}+B_{1}'+A_{2}'$ (Theorem 4.13)\\
\indent \ \ \ \ \ \ \ \ \ \ \ iff $A_{1}'+B_{1}+B_{2}'\equiv
B_{1}+B_{1}'+A_{2}'$ (by (*))\\
\indent \ \ \ \ \ \ \ \ \ \ \ iff $\ \ \ \ \ \ \
A_{1}'+B_{2}'\equiv
B_{1}'+A_{2}'$ (Theorem 4.13)\\
\end{rmk}

We then have;\\

\begin{theorem}{Transitivity of Linear Equivalence}\\

Let $C'$ be an algebraic curve. If $A,B,C$ are generalised
weighted sets on $C'$ of the same total multiplicity, then, if
$A\equiv B$ and $B\equiv C$, we must have that $A\equiv C$.

\end{theorem}

\begin{proof}
Let $\{A_{1},A_{2}\}$, $\{B_{1},B_{2}\}$ and $\{C_{1},C_{2}\}$ be
the effective weighted sets as given by Remarks 4.23. Then, by
Definition 4.24, we have that;\\

$(A_{1}+B_{2})\equiv (B_{1}+A_{2})$ and $(B_{1}+C_{2})\equiv
(C_{1}+B_{2})$\\

By Theorem 4.13, we have that;\\

$(A_{1}+B_{1}+B_{2}+C_{2})\equiv (C_{1}+B_{1}+B_{2}+A_{2})$\\

It then follows, by Definition 4.6, that there exists a
$g_{n}^{1}$, containing $(A_{1}+B_{1}+B_{2}+C_{2})$ and
$(C_{1}+B_{1}+B_{2}+A_{2})$ as weighted sets. Clearly
$(B_{1}+B_{2})$ is contained in the fixed branch contribution of
this $g_{n}^{1}$. Removing this fixed branch contribution, we
obtain;\\

$A_{1}+C_{2}\equiv C_{1}+A_{2}$\\

By Definition 4.24, we then have that $A\equiv C$ as required.

\end{proof}

It follows immediately from Theorem 4.12 and Theorem 4.26 that;\\

\begin{theorem}
Let $C$ be a projective algebraic curve, then $\equiv$ is an
equivalence relation on generalised weighted sets for $C$ of a
given total multiplicity.
\end{theorem}

\begin{rmk}
Again, the definition of linear equivalence that we have given for
generalised weighted sets on a smooth projective algebraic curve
$C$ is equivalent to the modern definition for divisors. More
precisely, we have;\\

Modern Definition; Let $A$ and $B$ be divisors on a smooth
projective algebraic curve $C$, then $A\equiv B$ iff $A-B=div(g)$,
for some $g\in L(C)^{*}.$\\

See, for example, p161 of \cite{Shaf} for relevant definitions and
notation. In order to show that our definition is the same, use
Remarks 4.8 and the following simple argument;\\

$A\equiv B$ as generalised weighted sets iff $A_{1}+B_{2}\equiv
B_{1}+A_{2}$\\

where $\{A_{1},A_{2},B_{1},B_{2}\}$ are the effective weighted
sets given by Definition 4.24. Then;\\

$A_{1}+B_{2}\equiv B_{1}+A_{2}$ iff
$(A_{1}+B_{2})-(B_{1}+A_{2})=div(g)$ $(g\in L(C)^{*})$\\

by Remarks 4.8, where $div(g)$ is the modern definition. By a
straightforward calculation, we have that;\\

$(A_{1}+B_{2})-(B_{1}+A_{2})=A-B$ as divisors or generalised
weighted\\
\indent \ \ \ \ \ \ \ \ \ \ \ \ \ \ \ \ \ \ \ \ \ \ \ \ \ \ \ \ \
\ \ \ \ \ \ \ \ \ \ \  sets.

Hence, the notions of equivalence coincide.

\end{rmk}

We also have;\\

\begin{theorem}{Linear Equivalence Preserved by Addition}\\

Let $C'$ be a projective algebraic curve and suppose that
$\{A,B,C,D\}$ are generalised weighted sets on $C'$ with;\\

$A\equiv B$ and $C\equiv D$\\

then;\\

$A+C\equiv B+D$\\

\end{theorem}

\begin{proof}
Let $\{A_{1},A_{2}\}$, $\{B_{1},B_{2}\}$, $\{C_{1},C_{2}\}$ and
$\{D_{1},D_{2}\}$ be effective weighted sets as given by Remarks
4.23 Then, by Definition 4.24, we have that;\\

$A_{1}+B_{2}\equiv B_{1}+A_{2}$ and $C_{1}+D_{2}\equiv
D_{1}+C_{2}$\\

Hence, by Theorem 4.13;\\

$A_{1}+B_{2}+C_{1}+D_{2}\equiv B_{1}+A_{2}+D_{1}+C_{2}$ $(*)$\\

We clearly have that;\\

$A+C=(A_{1}+C_{1})-(A_{2}+C_{2})$ and $B+D=(B_{1}+D_{1})-(B_{2}+D_{2})$\\

as an identity of generalised weighted sets. Moreover, as\\
$(A_{1}+C_{1}),(A_{2}+C_{2}),(B_{1}+D_{1})$ and $(B_{2}+D_{2})$
are all effective, we can apply Remarks 4.25 and $(*)$ to obtain
the result.
\end{proof}

We now make the following definition;\\

\begin{defn}
Let $G$ be a generalised weighted set on a projective algebraic
curve $C$, then we define $|G|$ to be the collection of
generalised weighted sets $G'$ with $G'\equiv G$. We define
$order(|G|)$ to be the total multiplicity (possibly negative) of
any generalised weighted set in $|G|$.
\end{defn}

\begin{rmk}
If $G$ is an \emph{effective} weighted set, the collection defined
by Definition 4.30 is \emph{not} the same as the collection given
by Definition 4.17, as it includes virtual weighted sets. Unless
otherwise stated, we will use Definition 4.17 for \emph{effective}
weighted sets. This convention is in accordance with the Italian
terminology.

\end{rmk}

We now show that the notions of linear equivalence introduced
in this section are birationally invariant;\\

\begin{theorem}
Let $\Phi:C_{1}\leftrightsquigarrow C_{2}$ be a birational map.
Let $A$ and $B$ be generalised weighted sets on $C_{2}$, with
corresponding generalised weighted sets $[\Phi]^{*}A$ and
$[\Phi]^{*}B$ on $C_{1}$. Then $A\equiv B$, in the sense of either
Definition 4.6 or 4.24, iff $[\Phi]^{*}A\equiv [\Phi]^{*}B$.

\end{theorem}

\begin{proof}
Suppose that $A\equiv B$ in the sense of Definition 4.6. Then,
there exists a $g_{n}^{r}$ on $C_{2}$ containing $A$ and $B$ as
weighted sets. By Theorem 3.14, there exists a corresponding
$g_{n}^{r}$ on $C_{1}$, containing $[\Phi]^{*}A$ and $[\Phi]^{*}B$
as weighted sets. Hence, again by Definition 4.6,
$[\Phi]^{*}A\equiv [\Phi]^{*}B$. The converse is similar, using
$[\Phi^{-1}]^{*}$. If $A\equiv B$ in the sense of Definition 4.24,
then the same argument works.

\end{proof}

As a result of this theorem, we introduce the following
definition;\\

\begin{defn}
Let $\Phi:C_{1}\leftrightsquigarrow C_{2}$ be a birational map.
Then, given a generalised weighted set $A$ on $C_{2}$, we
define;\\

$[\Phi]^{*}|A|=|[\Phi]^{*}A|$\\

where, in the case that $A$ is effective, $|A|$ can be taken
either in the sense of Definition 4.17 or Definition 4.30.

\end{defn}

\begin{rmk}
The definition depends only on the complete series $|A|$, rather
than its particular representative $A$. This follows immediately
from Definition 4.17, Definition 4.30 and Theorem 4.32.

\end{rmk}

We finally introduce the following definition;\\

\begin{defn}{Summation of Complete Series}\\

Let $A$ and $B$ be generalised weighted sets, defining complete
series $|A|$ and $|B|$, in the sense of Definition 4.30. Then, we
define the sum;\\

$|A|+|B|$\\

to be the complete series, in the sense of Definition 4.30,
containing all generalised weighted sets of the form $A'+B'$ with
$A'\in |A|$ and $B'\in |B|$. If $A$ and $B$ are \emph{effective}
weighted sets with $|A|$, $|B|$ taken in the sense of Definition
4.17, then we make the same definition for the sum in the sense of
Definition 4.17.\\

\end{defn}

\begin{rmk}
This is a good definition by Theorem 4.13 and Theorem 4.29.
\end{rmk}

\begin{defn}{Difference of Complete Series}\\

Let $A$ and $B$ be generalised weighted sets, defining complete
series $|A|$ and $|B|$, in the sense of Definition 4.30. Then, we
define the difference;\\

$|A|-|B|$\\

to be the complete series, in the sense of Definition 4.30,
containing all generalised weighted sets of the form $A'-B'$ with
$A'\in |A|$ and $B'\in |B|$. If $A$ and $B$ are \emph{effective}
weighted sets with $|A|$, $|B|$ taken in the sense of Definition
4.17, then we can in certain cases define a difference in the
sense of Definition 4.17. (This is called the residual series, the
reader can look at \cite{Sev} for more details)

\end{defn}

\begin{rmk}
This is again a good definition, for generalised weighted sets
$\{A,B\}$, it follows trivially from the previous definition and
the fact that $\{A,-B\}$ are also generalised weighted sets.
\end{rmk}

\end{section}


\begin{thebibliography}{99}

\bibitem{Aby} S.S. Abhyankar, Algebraic Geometry for Scientists
and Engineers, AMS Mathematical Surveys 35, (1990)\\

\bibitem{Hart} R. Hartshorne, Algebraic Geometry, Springer (1977)\\

\bibitem{Mum} D. Mumford. Red Book of Varieties and Schemes, Springer (1999)\\

\bibitem{Z} K. Peterzil and B. Zilber, Lecture Notes on Zariski Structures (1996)\\

\bibitem{depiro3} T. de Piro, A Non-Standard Bezout Theorem, AG/LO ArXiv (0406176), (2004).\\

\bibitem{depiro1} T. de Piro, A Theory of Branches for Algebraic Curves, MODNET preprint server, (2006).\\

\bibitem{depiro4} T. de Piro, Infinitesimals in a Recursively Enumerable Prime Model, LO ArXiv (0510412), (2005).\\

\bibitem{depiro2} T. de Piro, Zariski Structures and Algebraic Geometry, AG arXiv, math.AG/0402301, (2004).\\

\bibitem{Sev} F. Severi, Trattato di geometria algebrica, 1: Geometria delle serie lineari, Zanichelli, Bologna, (1926).\\

\bibitem{Shaf} I. Shafarevich, Basic Algebraic Geometry 1, Springer (1977).\\










\end{thebibliography}
\end{document}